\documentclass[11pt]{elsarticle}

\makeatletter
\def\ps@pprintTitle{%
 \let\@oddhead\@empty
 \let\@evenhead\@empty
 \def\@oddfoot{\centerline{\thepage}}%
 \let\@evenfoot\@oddfoot}
\makeatother

\usepackage{microtype}
\usepackage[]{algorithm2e}
\usepackage{chngpage}
\usepackage[english]{babel}
\usepackage[T1]{fontenc}
\usepackage[utf8]{inputenc}
\usepackage{amsfonts}
\usepackage{amsmath}
\usepackage{amssymb}
\usepackage{mathtools}
\usepackage[binary-units=true]{siunitx}
\usepackage{breqn}
\usepackage{stmaryrd}
\usepackage{sansmath}
\usepackage{amsthm}
\usepackage{bm}

\usepackage{caption}
\usepackage{subcaption}
\usepackage{diagbox}

\usepackage{rotating}

\usepackage{geometry}
\newgeometry{vmargin={25mm}, hmargin={25mm,25mm}}

\usepackage{multicol}
\usepackage{booktabs}
\usepackage{lscape}
\makeatletter
\newcommand\tabfill[1]{%
\dimen@\linewidth%
\advance\dimen@\@totalleftmargin%
\advance\dimen@-\dimen\@curtab%
\parbox[t]\dimen@{#1\ifhmode\strut\fi}%
}

\usepackage{multirow}

\usepackage{dcolumn}
\newcolumntype{d}[1]{D{.}{.}{#1}}

\usepackage{tikz}
\usepackage[]{graphicx}

\usetikzlibrary{shapes.geometric}

\renewcommand{\@algocf@capt@plain}{above}% formerly {bottom}

\definecolor{clr1}{rgb}{0.031, 0.270, 0.580}
\definecolor{clr2}{rgb}{0.129, 0.443, 0.709}
\definecolor{clr3}{rgb}{0.258, 0.572, 0.776}
\definecolor{clr4}{rgb}{0.419, 0.682, 0.839}
\definecolor{clr5}{rgb}{0.619, 0.792, 0.882}

\usepackage[unicode=true,
            bookmarks=true,
            bookmarksnumbered=false,
            bookmarksopen=false,
            breaklinks=true,
            pdfborder={0 0 1},
            backref=false,
            colorlinks=true,
            linkcolor=black,
            urlcolor=theblue,
            citecolor=theblue,
            hyperfootnotes=false]{hyperref}
\usepackage{lineno}
\hypersetup{pdftitle={},
            pdfauthor={The Authors}}

\let\oldref\ref
\renewcommand{\ref}[1]{(\oldref{#1})}

\allowdisplaybreaks

\usepackage{placeins}
\begin{document}

\newcommand{\hwsolid}{\raisebox{2pt}{\tikz{\draw[black,solid,line width=0.9pt](0,0) -- (5mm,0);}}}
\newcommand{\hwdashdotted}{\raisebox{2pt}{\tikz{\draw[black,dashed,line width=0.9pt](0,0) -- (5mm,0);}}}
\newcommand{\hwdotted}{\raisebox{2pt}{\tikz{\draw[black,dotted,line width=0.9pt](0,0) -- (5mm,0);}}}

\newcommand{\hwcircle}{\raisebox{0pt}{\tikz{\node[circle,draw,fill=none,text width=2mm, inner sep=0pt] at (0,0){};}}}
\newcommand{\hwtriangle}{\raisebox{0pt}{\tikz{\node[isosceles triangle,draw,fill=none,rotate=90,text width=2mm, inner sep=0pt] at (0,0){};}}}
\newcommand{\hwsquare}{\raisebox{0pt}{\tikz{\node[regular polygon,regular polygon sides=4,draw,fill=none,text width=1.25mm, inner sep=0pt] at (0,0){};}}}

\newcommand{\hwcirclefilled}{\raisebox{0pt}{\tikz{\node[circle,draw,fill=black,text width=2mm, inner sep=0pt] at (0,0){};}}}
\newcommand{\hwtrianglefilled}{\raisebox{0pt}{\tikz{\node[isosceles triangle,draw,fill=black,rotate=90,text width=2mm, inner sep=0pt] at (0,0){};}}}
\newcommand{\hwsquarefilled}{\raisebox{0pt}{\tikz{\node[regular polygon,regular polygon sides=4,draw,fill=black,text width=1.25mm, inner sep=0pt] at (0,0){};}}}

% \graphicspath{{figures/}} 

\title{Hybridized Implicit-Explicit Flux Reconstruction Methods}

\begin{frontmatter}
\author[one]{Carlos A. Pereira\corref{mycorrespondingauthor}}
\cortext[mycorrespondingauthor]{Corresponding author}
\ead{carlos.pereira@concordia.ca}

\author[one]{Brian C. Vermeire}
\ead{brian.vermeire@concordia.ca}

\address[one]{Department of Mechanical, Industrial, and Aerospace Engineering \\Concordia University \\Montreal, QC. Canada}

\begin{abstract}
For turbulent problems of industrial scale, computational cost may become prohibitive due to the stability constraints associated with explicit time discretization of the underlying conservation laws. On the other hand, implicit methods allow for larger time-step sizes but require exorbitant computational resources. Implicit-explicit (IMEX) formulations combine both temporal approaches, using an explicit method in nonstiff portions of the domain and implicit in stiff portions. While these methods can be shown to be orders of magnitude faster than typical explicit discretizations, they are still limited by their implicit discretization in terms of cost. Hybridization reduces the scaling of these systems to an effective lower dimension, which allows the system to be solved at significant speedup factors compared to standard implicit methods. This work proposes an IMEX scheme that combines hybridized and standard flux reconstriction (FR) methods to tackle geometry-induced stiffness. By using the so-called transmission conditions, an overall conservative formulation can be obtained after combining both explicit FR and hybridized implicit FR methods. We verify and apply our approach to a series of numerical examples, including a multi-element airfoil at Reynolds number $1.7$ million. Results demonstrate speedup factors of four against standard IMEX formulations and at least 15 against standard explicit formulations for the same problem.
\end{abstract}
\end{frontmatter}
\section{Introduction}

The behaviour of turbulent flows at industrial scales is inherently unsteady and complex. Hence, appropriate temporal methods must be chosen to advance the solution in time. These methods are generally classified as explicit, implicit, and implicit-explicit (IMEX). Explicit methods, such as explicit Runge-Kutta (ERK) schemes, calculate the state of the system at a later time or stage from a known value at the current time or stage. They are generally easy to implement, do not require significant memory, and have a short computation time per time step~\cite{hairerSolvingOrdinaryDifferential1993,butcherPracticalRungeKutta2009}. However, the maximum allowable time-step size is determined by the relationship between the spatial discretization and a constrained stability region that dictates the maximum stable Courants-Fredrichs-Lewys (CFL) number. Hence, they are generally suitable for nonstiff problems. On the other hand, implicit methods calculate the future state of the system via the solution of a coupled system of equations. Unlike explicit methods, implicit schemes can be unconditionally linearly stable~\cite{hairerSolvingOrdinaryDifferential1996}, and the value of the time-step size is typically chosen based on accuracy rather than stability. In general industrial applications, the solution of an implicit method must be obtained via coupled nonlinear solvers.

With advances in high-order methods, computations of more complex problems are becoming more feasible. Typical industrial applications of flow over aircraft wings are in the range of seven/eight-figure Reynolds numbers. The higher the value of this parameter, the denser the computational grid must be to capture the resulting thin boundary layer and smallest turbulent length scales. This poses an important challenge for the use of explicit methods due to their stability constraints, as the largest allowable time-step size becomes prohibitively small for these problems. In implicit methods, however, we are required to perform linearization, storage, preconditioning, and solution of large nonlinear systems, which scale very rapidly with the order of accuracy at $O(p^d)$, where $p$ is the polynomial degree representing the solution, and $d$ is the dimension of the problem. For typical problems of engineering interest, wall-resolved LES simulations using purely implicit methods necessitate infeasible computational resources.

IMEX methods combine the use of explicit and implicit methods to balance efficiency, stability, and accuracy, using the best features of each approach~\cite{ascherImplicitexplicitRungeKuttaMethods1997,shoeybi2010adaptive,hundsdorfer2003implicit}. Generally, they are designed such that an implicit method is used to solve the stiff components of the underlying equations and an explicit method is used on the nonstiff portion. Then, the two are paired to enable time integration that maintains conservation and a desired order of accuracy. With IMEX methods, the maximum CFL condition can be relaxed in the explicit part, and the simulation can be run more efficiently. IMEX methods were initially devised to time-split the convection-diffusion-reaction operator into stiff and nonstiff components~\cite{ascherImplicitexplicitRungeKuttaMethods1997}. More recently, they have been applied to tackle so-called geometry-induced stiffness~\cite{kanevskyApplicationImplicitExplicit2007,vermeire2015adaptive}. This occurs when there exist largely disparate cell sizes in a single computational domain, as in, for instance, the cells within the boundary layer of an airfoil at high-Reynolds numbers compared to those in the far-field~\cite{kanevskyApplicationImplicitExplicit2007}. IMEX methods have also been developed for specific applications such as in pseudo-time stepping for incompressible flows~\cite{vermeire2020optimal} and combined with $p$-multigrid methods~\cite{jameson1983solution,fidkowski2005p} for convergence acceleration. IMEX schemes have been shown to provide speedups over explicit methods by one to two orders of magnitude, thereby reducing the computational cost of LES simulations~\cite{vermeire2015adaptive}. 

More recently, novel optimal explicit and implicit-explicit methods were introduced. These are obtained by reshaping stability polynomials of typical explicit Runge-Kutta (RK) methods to increase their maximum allowable time-step size for specific spatial discretization methods and conservation laws~\cite{ketchesonOptimalStabilityPolynomials2012,kubatkoOptimalStrongStabilityPreservingRunge2014,verm-ac}.   Optimal explicit RK schemes for the FR methods have been devised and analyzed in~\cite{pereiraFullyDiscreteAnalysisHighOrder2020}, and optimal IMEX RK schemes have been formulated in~\cite{vermeire2021accelerated}, known as accelerated IMEX (AIMEX) methods. These optimal IMEX formulations have been shown to yield speedup factors in excess of two over standard IMEX methods. While both IMEX and AIMEX methods have shown significant speedups over conventional explicit time-stepping, it is known that these schemes usually spend the majority of their computation time solving the nonlinear systems that result from the implicit component~\cite{vermeire2015adaptive}.

Hybridization can reduce the cost of the implicit solver, which typically takes most of the computation time in implicit and IMEX formulations. It was introduced in the context of DG methods by Cockburn et al.~\cite{cockburnUnifiedHybridizationDiscontinuous2009}, associated with the static condensation of finite-element methods of de Veubeke~\cite{de2001displacement} and extended to the family of FR schemes in~\cite{pereiraPerformanceAccuracyHybridized2022}, also known as hybridized Flux Reconstruction (HFR) schemes. HFR methods define a new unknown on the faces of the elements, known as the trace variable, which effectively decouples interelement information in the Jacobian matrix. Then, with transmission conditions, the problem is globally defined in terms of the trace variable and later reduced via static condensation~\cite{cockburnStaticCondensationHybridization2016}. This leads to a Jacobian matrix, the size of which grows proportionally to $O(p^{d-1})$ instead of the typical $O(p^{d})$ of standard implicit schemes. 

In standard HFR methods, the trace unknown can be chosen to be discontinuous. Subsets of these methods have been introduced by modifying the functional space of the trace variable. For example, by enforcing continuity on the skeleton of the domain, the number of globally coupled degrees of freedom decreases. This is known as the embedded discontinuous Galerkin (EDG) method initially introduced in the context of linear shell problems~\cite{guzeyEmbeddedDiscontinuousGalerkin2007} and embedded FR methods (EFR) in~\cite{pereiraPerformanceAccuracyHybridized2022}. An analysis of this method was performed in~\cite{cockburnAnalysisEmbeddedDiscontinuous2009}, where it was shown that although embedded methods result in smaller linear systems, they are only conservative on the dual volumes~\cite{kamenetskiyRelationEmbeddedDiscontinuous2016}. 

Combining hybridization with standard formulations has shown to provide significant benefits to reducing stiffness in partial differential equations via IMEX formulations. Previous works on IMEX methods with hybridization have been developed to tackle stiffness associated with shallow water systems~\cite{kangIMEXHDGDGCoupled2020}. This was done to separate the faster gravity waves from the nonlinear advection operator. However, operator-splitting applications typically struggle to compete with purely explicit methods. Coupling hybridized and standard methods to tackle stiffness associated with geometry has not been explored before. The objective in this work is to introduce a conservative time accurate pairing of hybridzed and standard flux reconstruction formulations for problems with geometry-induced stiffness.

This paper is structured as follows. Section 2 introduces the flux reconstruction approach, in both its standard and hybridized formulations. Section 3 describes the global system resulting from the hybridized formulation and the details of the time discretization. In section 4, implicit-explicit coupling of standard and hybridized methods is introduced by employing transmission conditions along the IMEX interface. Section 5 presents verification of the proposed approach as well as performance and validation studies via a series of numerical examples including flow over a multi-element airfoil at $\operatorname{Re}=1.7\times10^6$. Finally, conclusions and future work are presented in Section 6.

\section{The Flux Reconstruction Method}\label{sec:hfr}

Consider the following conservation law
\begin{equation}
 \frac{\partial u}{\partial t} + \nabla \cdot \bm F(u) = 0~~\text{in}~\Omega,
\end{equation}
where $\Omega$ is a bounded subset of $\mathbb{R}^d$ with boundary $\partial\Omega\in \mathbb{R}^{d-1}$ and $d$ dimensions, $u$ is the conserved quantity, $\bm F=\bm F(u)$ is the flux, and $t$ is time. 

Define the computational domain by a partition of nonoverlapping, conforming elements $\Omega_k$ such that  $\mathcal{T}_h = \{\Omega_k\}$. The boundary of each element is defined by $\partial\Omega_k=\{f\}$ with $|\partial\Omega_k|=N_f$ which we collectively include in the set $\partial \mathcal{T}_h=\{\partial \Omega_k: \Omega_k \in \mathcal{T}_h\}$. In this set, the two faces belonging to neighbouring elements coexist. The intersection of all element faces in the domain defines the skeleton of the grid $\varepsilon_h=\varepsilon^\partial\cup\varepsilon^0=\{\bar f\}$, where $\varepsilon^\partial$ refers to the boundary faces, and $\varepsilon^0$ represents the interior faces. A relationship between an element's face and its global position in the computational grid is given by $\bar f\in \varepsilon_h,~\bar f = (f\in\partial\Omega_k) \cap \varepsilon_h$.

The relationship between the reference and physical space is given via invertible one-to-one mapping functions $\bm {\mathcal{M}}_k( \tilde{\bm x})$ for each element. Hence, any physical coordinate $\bm x$ in $\Omega_k$ can be obtained from a reference location via
\begin{equation}
  \bm x = \bm {\mathcal{M}}_k( \tilde{\bm x}) = \sum_{i=1}^{N_g} M_i( \tilde{\bm x}) \bm x^g_i,
\end{equation}
where $M_i$ is a mapping function defined by $N_g$ mapping points $\{\bm x^g_i\}$. The Jacobian matrix of these transformations is defined by $\bm J_k(\tilde{\bm x})$ and its determinant by $J_k(\tilde{\bm x})$. The relationship between physical and reference quantities in the conservation law can be written for a time-invariant formulation~\cite{zwanenburgEquivalenceEnergyStable2016}
\begin{align}
  \tilde{u}^h_k &= \tilde{u}^h_k(\tilde{\bm  x}, t) = J_k u^h_k(\bm{\mathcal{M}}_k(\tilde{\bm  x}), t),\\
  \tilde{\bm F}^h_k &= \tilde{\bm F}^h_k(\tilde{\bm  x}, t) = J_k \bm J_k^{-1} \bm F^h_k(\bm{\mathcal{M}}_k(\tilde{\bm  x}), t),
 \label{eq:d9sfxda}
\end{align}
which allows us to evaluate the evolution of the solution in physical space via
\begin{align}
  \frac{\partial u^h_k}{\partial t} + \frac{1}{J_k}\tilde \nabla \cdot\tilde{\bm F}^h_k = 0,
  \label{eq:3289u4}
\end{align}
where $\tilde\nabla$ is the reference space divergence operator. Within each element, $N_s$ interior solution points $\{\tilde{\bm x}_i^s\}_{i=1}^{N_s}$ define globally discontinuous polynomials of degree $p$. Hence, we can represent polynomials for conserved variables using nodal basis functions $\varphi(\tilde{\bm x})$
\begin{equation}
  u^h_k (\tilde{\bm x},t) = \sum_{i=1}^{N_s} U_{k,i}(t) \varphi_i(\tilde{\bm x}).
  \label{eq:frsolution}
\end{equation}
A discontinuous flux function can also be represented using nodal values, leading to a polynomial in the same space as the solution, that is
\begin{equation}
  \tilde{\bm F}^{hD}_k (\tilde{\bm x},t) = \sum_{i=1}^{N_s} \tilde{\bm F}_{k,i}(t) \varphi_i(\tilde{\bm x}),
  \label{eq:frfluxd}
\end{equation}
where $\tilde{\bm F}_{k,i}(t)$ is the transformed flux evaluated at the solution point $\tilde{\bm x}_i^s$. For conservation, a globally $C_0$-continuous flux is determined by adding the following term
\begin{equation}
\tilde {\bm F}_{k}^{hC}(\tilde{\bm x}, t) = \sum_{f=1}^{N_f}\sum_{m=1}^{N_{r,f}} \bm g_f^m(\tilde{\bm x}) \left[{\tilde{ H}}(\tilde{\bm x})_{k,f}\right]_{\tilde{\bm x} = \tilde{\bm x}_{f,m}^r}.
\label{eq:fluxcorrection}
\end{equation}
Several new variables are introduced in this context, including the correction functions. Among these are the correction functions $\bm g$, which belong to a Raviart-Thomas space and are  defined on each of the flux points $\{\tilde{\bm x}_{f,i}^r\}_{i=1}^{N_{r,f}}$ at each face $f$. These functions satisfy
\begin{equation}
  \tilde{\bm n}_f^m \cdot \bm{g}^n_l(\tilde{\bm x}_f^m) = \delta_{fl}\delta_{mn}.
  \label{eq:dsklnf}
\end{equation}
In addition, we have introduced the flux interface jump, defined as
\begin{equation}
  \tilde{H}_{k,f}(\tilde{\bm x}) = \tilde{\hat{\bm{\mathfrak{F}}}}_{k,f}\cdot \tilde{\bm n}_f  - \tilde{\bm F}_{k,f}^{hD}\cdot\tilde{\bm n}_f,
\end{equation}
where $\tilde{\bm F}_{k,f}^{hD}$ is the transformed discontinous flux evaluated at the interface, $\tilde{\bm n}_f$ is the reference outward normal vector function and $\tilde{\hat{\bm{\mathfrak{F}}}}_{k,f}$ is the Riemann flux at the interface. The relationship between physical and reference space for the common flux is~\cite{zwanenburgEquivalenceEnergyStable2016}
\begin{equation}
  \tilde{\hat{\bm{\mathfrak{F}}}}_{k,f}\cdot \tilde{\bm n}_f = J_{k,f} \hat{\bm{\mathfrak{F}}}_{k,f}\cdot {\bm n}_{k,f},
  \label{eq:8923i}
\end{equation}
where $J_{k,f}$ is the Jacobian determinant of the face. There remains to be defined both the corrected solution and flux interface variables, which lead to different types of discretizations. Depending on the choice of the common solution and fluxes, the numerical method will have different stability properties, and the stencil will also vary. Hence, the corrected flux becomes
\begin{equation}
  \tilde{\bm F}_k^h = \tilde{\bm F}_k^{hD} + \tilde{\bm F}_k^{hC}.
  \label{eq:correctedflux}
\end{equation}
\subsection{Standard FR formulation}~\label{sec:fr}
In standard formulations, the numerical flux at the interface typically takes the following form for the inviscid or convective component
\begin{equation}
    \hat{\bm{\mathfrak{F}}} (u_-, u_+) \cdot \bm n_- =  \frac{1}{2}\left[\bm F(u_-) + \bm F(u_+)\right]\cdot \bm n_- + \frac{s^{\text{FR}} }{2}(u_- - u_+),
\end{equation}
where $s^{\text{FR}}$ is a stabilization parameter in the standard formulation. This parameter is chosen depending on the physics of the problem. Typical values (or matrices in systems of equations) lead to the Rusanov or local Lax-Friedrichs (LLF)~\cite{rusanovCalculationInteractionNonsteady1962,laxWeakSolutionsNonlinear1954}, Roe~\cite{roeApproximateRiemannSolvers1981}, and HLL methods~\cite{toro1994restoration}. Throughout this work, we employ the Rusanov/LLF fluxes with $s^{\text{FR}} = \bar \lambda^{\text{FR}}$, with $\bar\lambda^{\text{FR}}$ the maximum local wave speed in the system. 

After choosing a suitable interface formulation, we apply the divergence on the corrected flux in Equation~\eqref{eq:correctedflux} and sum over all elements, which yields the following system for the standard FR method
\begin{subequations}\begin{align}
  \sum_{\Omega_k\in \mathcal{T}^h}\frac{\partial u^h_k}{\partial t} + \frac{1}{J_k}\sum_{i=1}^{N_s} \tilde{\bm F}_{k,i} \cdot \tilde \nabla \varphi_i(\tilde{\bm x}) +  \frac{1}{J_k}\sum_{f=1}^{N_f}\sum_{m=1}^{N_{r,f}} \tilde\nabla\cdot \bm g_f^m(\tilde{\bm x}) \left[{\tilde{H}}(\tilde{\bm x})_{k,f}\right]_{\tilde{\bm x} = \tilde{\bm x}^r_{f,m}}&= 0\label{eq:freq1}.
\end{align}\label{eq:freq}\end{subequations}
Since the definitions of the interface variables are explicitly defined as a function of left and right states, this system can directly be solved both implicitly and explicitly. Similar to this methodology, we now demonstrate the steps for the hybridized method. 

\subsection{Hybridized FR formulation}

In addition to the conserved variable, we introduce an approximation to $u^h$ on the skeleton of the computational grid such that at any face $\bar f\in\varepsilon^h$, a degree-$p$  polynomial can be obtained via
\begin{equation}
  \hat{u}^h_{\bar f} (\tilde{\bm x},t) = \sum_{i=1}^{N_{r,f}}  {\hat{U}}_{\bar f,i}(t) \phi_i(\tilde{x}),
  \label{eq:rwepsdf}
\end{equation}
which is the so-called trace variable. Here, $N_{r,f}$ is the number of flux points in face $f$, which we consider equal to the number of trace points at a given face, and $\phi$ is a $d-1$-dimensional trace basis function. Hybridization is then achieved by considering the following form of the common fluxes
\begin{equation}
 \hat{{\bm{\mathfrak{F}}}}_{k,f} = \bm F(\hat{u}^h_{\bar f}) + s(u^h_{k,f} - \hat{u}^h_{\bar f}) {{\bm n}}_{k,f},
 \label{eq:hfrriemann}
\end{equation}
with ${{\bm n}}_{k,f}$ the physical outward unit normal vector and $s$ a stabilization parameter~\cite{nguyenImplicitHighorderHybridizable2009,pereiraPerformanceAccuracyHybridized2022}, which is chosen following the physics of the conservation law. These fluxes depend on information from a single element and the corresponding trace. Contrary to the standard FR implementation, we do not strongly enforce conservation at the interface since $\hat{\bm{\mathfrak F}}_-$ need not be equal to $\hat{\bm{\mathfrak F}}_+$ in the general case. Hence, we seek to enforce an additional statement satisfying discrete global conservation, i.e., 
\begin{equation}
  \sum_{\bar f\in \varepsilon^h_0} \int_{\bar f} \llbracket \hat{\bm{\mathfrak{F}}}\rrbracket_{\bar f} \phi ds+ \sum_{\bar f\in \varepsilon^h_\partial}  \int_{\bar f}  {\mathfrak{F}}^{\text{BC}}_{\bar f} \phi ds = 0,
  \label{eq:5treerf}
\end{equation}
and provides closure to the system. In these equations, we have separated the interior and normal boundary fluxes, the latter of which we denote ${\mathfrak{F}}^{\text{BC}}_{\bar f}$. The jump operator is defined at an interface $\bar f$ between elements $\Omega_+$ and $\Omega_-$ as
\begin{equation}
  \llbracket\bm f \rrbracket_{\bar f} = \bm f_{k^+,f^+} \cdot { {\bm{n}}}_{k^+,f^+} + \bm f _{k^-,f^-} \cdot {\bm{n}}_{k^-,f^-}.
  \label{eq:53904jn}
\end{equation}
After summing over all elements, we can state the hybridized form of the flux reconstruction approach as follows
\begin{subequations}\begin{align}
  %\sum_{\Omega_k\in \mathcal{T}^h}\tilde{\bm q}_k^h - \sum_{i=1}^{N_s}  U_{k,i} \tilde\nabla \varphi_i(\tilde{\bm x}) - \sum_{f=1}^{N_f} \sum_{m=1}^{N_{r,f}} \tilde{\bm n}_f^m\cdot \tilde \nabla \cdot \bm g_f^m(\tilde{\bm x})\left[\mathfrak{U}_{k,f} - u_{k,f}^h\right]_{\tilde{\bm x} = \tilde{\bm x}^r_{f,m}} &= 0\label{eq:hfr0},\\
  \sum_{\Omega_k\in \mathcal{T}^h}\frac{\partial {u}^h_k}{\partial t} +  \frac{1}{J_k}\sum_{i=1}^{N_s} \tilde{\bm F}_{k,i} \cdot \tilde \nabla \varphi_i(\tilde{\bm x}) +  \frac{1}{J_k}\sum_{f=1}^{N_f}\sum_{m=1}^{N_{r,f}} \tilde\nabla\cdot \bm g_f^m(\tilde{\bm x}) \left[{\tilde{H}}(\tilde{\bm x})_{k,f}\right]_{\tilde{\bm x} = \tilde{\bm x}^r_{f,m}}&= 0,\label{eq:hfr1}\\
 \sum_{\bar f\in \varepsilon^h_0} \int_{\bar f} \llbracket \hat{\bm{\mathfrak{F}}}\rrbracket_{\bar f} \phi ds+ \sum_{\bar f\in \varepsilon^h_\partial}  \int_{\bar f}  {\mathfrak{F}}^{\text{BC}}_{\bar f} \phi ds &= 0,\label{eq:hfr2}
\end{align}\label{eq:hfr}\end{subequations}
where we have readily taken the divergence of the flux and its correction to arrive at~\eqref{eq:hfr1}. Typically, hybridized methods make use of discontinuous or globally continuous function spaces for the trace variable, which can be respectively defined by
\begin{subequations}
\begin{align}
  \mathbb{M}^h_p &= \{\phi \in L_2(\varepsilon^h)~:~\phi|_{\bar f} \in \mathbb{P}^p(\bar f),~\forall \bar f \in \varepsilon^h\},\label{eq:hfrspaceshfr}\\
  \bar{\mathbb{M}}^h_p &= \mathbb{M}^h_p \cap C^0(\varepsilon^h).
\end{align}\label{eq:hfrspaces}\end{subequations}
The choice of nodal basis functions for the trace variable leads to different types of hybridizations. Here ${\mathfrak{F}}^{\text{BC}}_{\bar f}$ is the normal boundary flux. 
\begin{figure}[htb]
  \centering
  \begin{subfigure}[b]{0.49\textwidth}
    \includegraphics[width=\textwidth]{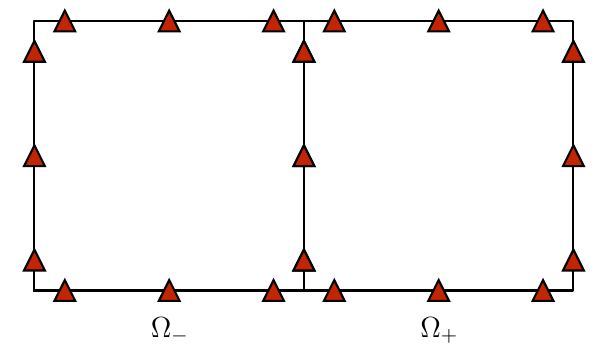}
    \caption{HFR}
  \end{subfigure}
  \begin{subfigure}[b]{0.49\textwidth}
    \includegraphics[width=\textwidth]{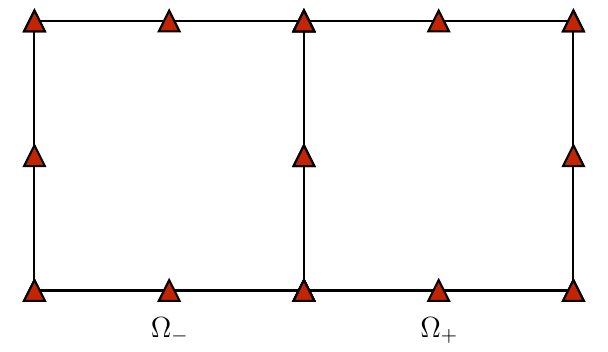}
    \caption{EFR}
  \end{subfigure}
  \caption[]{Trace variable location in an HFR (left) and EFR (right) discretization considering a $p=2$ scheme on the skeleton of two neighbouring quadrilateral elements}
  \label{fig:tracepoints}
\end{figure}
In this work, we consider two types of hybridization. First, we consider discontinuous trace nodal basis functions $\phi\in\mathbb{M}^h_p$. This leads to the HFR method. We also define globally continuous nodal basis functions on interior faces $\varepsilon^h_0$ $\phi\in\bar{\mathbb{M}}^h_p$ and discontinuous on boundary faces $\varepsilon_\partial^h$. This leads to the interior-embedded flux reconstruction (IEFR) scheme, which we will denote EFR for simplicity. A diagram of the resulting trace points is shown for two neighbouring elements $\Omega_-$ and $\Omega_+$ in Figure~\ref{fig:tracepoints} for a $p=2$ scheme. In this work, we choose the solution points inside the element to be those of the Gauss-Legendre quadrature. In the case of EFR, Gauss-Lobatto-Legendre (GLL) points are used at the faces. When different sets of points are used in the EFR method for the volume and faces, the correction functions are still generated using a tensor-product formulation, and the lifting operation includes an additional interpolation operator. See~\cite{pereiraPerformanceAccuracyHybridized2022}.

\section{The Global System}\label{sec:global}
In this section, we discuss the implicit temporal discretization and global system of the hybridized system. Later on in this manuscript, we will discuss how to couple this system with an explicit standard method to obtain a hybridized IMEX formulation.

\subsection{Temporal Discretization}
Solving a hybridized problem is generally done implicitly. We make use of second-order two-stage implicit SDIRK methods to advance the solution in time, which has the Butcher tableau
\begin{equation}
  \begin{tabular}{c|cccc}
  $\bm c$   & $\bar{\bm A}$ \\
  \hline
   & $\bm b$\\
  \end{tabular}
    =   
  \begin{tabular}{c|cccc}
   $1-\frac{\sqrt{2}}{2}$    & $1-\frac{\sqrt{2}}{2}$ \\
   $1$    & $\frac{\sqrt{2}}{2}$  & $1-\frac{\sqrt{2}}{2}$ \\
  \hline
           & $\frac{\sqrt{2}}{2}$  & $1-\frac{\sqrt{2}}{2}$\\
  \end{tabular}.
\end{equation}
Furthermore, after hybridization, we can write Equations~\eqref{eq:hfr}
\begin{subequations}\begin{align}
  \frac{\partial {\bm u}}{\partial t} + f ({\bm u}, \hat {\bm u}) &= 0, \label{eq:099io} \\
  g ({\bm u},\hat {\bm u}) &= 0,
  \label{eq:9mi90j}
\end{align}\end{subequations}
where $\bm u\in\mathbb{R}^N$, $\hat{\bm u}\in\mathbb{R}^{\bar N}$ are the interior and trace solution vectors, $f$ is the hybridized spatial discretization function associated with the first equation in~\eqref{eq:hfr1}, and $g$ is the residual associated with the transmission conditions~\eqref{eq:hfr2} . Hence, at the $i$-th stage, the solution is given by the system
\begin{subequations}
\begin{align}
  {\bm u}^{(i)} - {\bm u}^{(t)} - \sum_{j=1}^{i-1} a_{ij} \Delta t f(\bm u^{(j)},\hat {\bm u}^{(j)})  - a_{ii} \Delta t f({\bm u}^{(i)},\hat {\bm u}^{(i)}) =0&,\\
  g({\bm u}^{(i)},\hat {\bm u}^{(i)})  = 0&,
\end{align}
\end{subequations}
where information associated with the diagonal coefficients is unknown, and with the off-diagonal coefficients and the current time $t$ is known. These equations can be simplified to
\begin{subequations}\begin{align}
  {\bm u}^{(i)} - {\bm u}^{(t*)} - a_{ii} \Delta t f({\bm u}^{(i)},\hat {\bm u}^{(i)})=0&,\label{eq:9oind}\\
  g({\bm u}^{(i)},\hat {\bm u}^{(i)})  = 0&\label{eq:23nffsd},
\end{align}\label{eq:dklasd}\end{subequations}
where
\begin{equation}
  {\bm u}^{(t*)} = \bm u^{(t)} + \sum_{j=1}^{i-1} a_{ij} \Delta t f({\bm u}^{(j)},\hat {\bm u}^{(j)}),
\end{equation}
is the known information at this stage $i$.

\subsection{Block Formulation}
Hybridized unsteady problems can be written\begin{subequations}
\begin{align}
h(\bm u, \hat {\bm u}) = 0&,\\
g(\bm u, \hat {\bm u}) = 0&,
\end{align}
\end{subequations}
where $\bm u\in\mathbb{R}^{N\times N_s}$, $f$ is the hybridized spatial discretization function associated with Equation~\eqref{eq:hfr1}, and $g$ is the residual associated with the flux conservation statement in Equation~\eqref{eq:hfr2}~\cite{fernandezEntropystableHybridizedDiscontinuous2019a}. At the $n$-th Newton iteration, we have
\begin{align}
\begin{bmatrix}
 \bm A^n & \bm B^n \\ 
 \bm C^n & \bm D^n 
\end{bmatrix}
\begin{bmatrix}
\delta \bm  u^n \\ 
\delta \hat{\bm  u}^n
\end{bmatrix}
= 
\begin{bmatrix}
\bm r^n \\ \bm s^n
\end{bmatrix},
\end{align}
where $\delta \bm  u^n$ and $\delta \hat{\bm  u}^n$ refer to the update vector of internal and trace solution points at this iteration. Due to the discontinuous nature of the interior solution and decoupling neighbouring elements' interior solutions, $\bm A$ is block-diagonal. This proves efficient when we reduce the problem via static condensation and solve the condensed global problem
\begin{equation}
  \bm L^n \delta \hat{\bm u}^n = \bm t^n,
\end{equation}
where $\bm L = \bm D - \bm C \bm A^{-1} \bm B$ and $\bm t = \bm s - \bm C \bm A^{-1} \bm r$. Then, the solution can be obtained locally for each element from
\begin{equation}
  \delta \bm u_k = \bm A_k^{-1}(\bm r_k - \bm B_k \delta \hat {\bm u}_k).
  \label{eq:localproblems}
\end{equation}
Constructing the global operators can be done efficiently on a per-element basis, i.e.
\begin{align}
    \label{eq:assembly} \bm L_{i,j} &= \bm L_{i,j} + \bm L^k_{\bar i,\bar j}, \\
    \label{eq:assembly2} \bm t_{i} &= \bm t_{i} + \bm t^k_{\bar i},
\end{align}
with elemental matrices $\bm L^k$ and $\bm t^k$ thereby defined by
\begin{align}
   \bm L^k &:= \bm D^k - \bm C^k(\bm {A}^{-1})^k\bm B^k, \\
   \bm t^k &:= \bm s^k - \bm C^k (\bm A^{-1})^k\bm r^k,
\end{align}
and the indices $\bar i,~\bar j$ are associated with a surjective mapping of the element's flux points to the global trace points. The elemental blocks $\bm A^k,~\bm B^k$ can be specifically defined as follows
\begin{subequations}
\begin{align}
    \label{eq:elemA}&\bm A_{k,ij} = \delta_{ij} - \Delta t a_{ii} \frac{1}{J_{k,i}}\left[\sum_{g=1}^{N_s} \tilde \nabla \varphi_g(\tilde{\bm x}) \cdot  \frac{\partial\tilde{\bm F}_{k,g}}{\partial u_{k,j}}  + \sum_{f=1}^{N_f}\sum_{m=1}^{N_{r,f}} \tilde\nabla\cdot \bm g_f^m (\tilde{\bm x}) \frac{\partial \tilde H_{k,f}}{\partial u_{k,j}}(\tilde{\bm x}^r_{f,m})\right]_{\tilde{\bm x} = \tilde{\bm x}_i^s},\\ 
    &\bm B_{k,il} = \frac{\Delta t a_{ii}}{J_{k,i}}\left[\sum_{f=1}^{N_f}\sum_{m=1}^{N_{r,f}} \tilde\nabla\cdot \bm g_f^m (\tilde{\bm x}) \frac{\partial \tilde H_{k,f}}{\partial \hat u_{k,l}}(\tilde{\bm x}^r_{f,m})\right]_{\tilde{\bm x} = \tilde{\bm x}_i^s},
\end{align}\label{eq:elemsystem}\end{subequations}
and the remaining blocks can be further structured face-by-face within each element and then assembled into a single block. This way, we can easily accommodate faces with different numbers of trace points, which results from $p$-adaptive algorithms. That is, $\bm C^k$ and $\bm D^k$ can be trivially assembled from
\begin{subequations}
\begin{align}
  &\bm M_{k,f}^{-1} \bm C^{k,f}_{qi} =  \frac{\partial \tilde{\hat{\bm{\mathfrak{F}}}}_{k,f}}{\partial u_{k,j}} (\tilde {\bm x}^r_{f,q}) \cdot \tilde{\bm n}_f^q, \\
  &\bm M_{k,f}^{-1} \bm D^{k,f}_{qt} = \frac{\partial \tilde{\hat{\bm{\mathfrak{F}}}}_{k,f}}{\partial \hat u_{k,f,t}} (\tilde {\bm x}^r_{f,q}) \cdot \tilde{\bm n}_f^q,
\end{align}\label{eq:elemsystemB}\end{subequations}
where the indices $i=1,\ldots,N_s$, $l=1,\dots,N_r$ and $p,q=1,\ldots,N_{r,f}$, $N_r=\sum_{f}^{N_f} N_{r,f}$. Note that $\hat u_{k,l}$ refers to the trace living on the $l$-th point of element $\Omega_k$, and $\hat u_{k,f,t}$ is the trace at the $t$-th point of the $f$-th face in element $\Omega_k$. In addition, $\bm M_{k,f}$ contains the element local face mass matrices $\bm M_{k,f,ij}=\int_{\bar f}{\phi_i\phi_j}d\bar f$ as diagonal blocks blocks. Similarly, the vectors $\bm r^k,~\bm s^{k}$ evaluate the right-hand-side functions in the Newton algorithm.

\section{IMEX Formulation}

A largely disparate range of element sizes typically appears in simulations at high Reynolds numbers, with large elements in the far field and elements orders of magnitude smaller in proximity to walls. This introduces what is referred to as geometry-induced stiffness since the maximum time-step size for explicit simulations is dictated by the smallest element size. For high-order schemes, using explicit methods can be very restrictive in terms of the allowable time-step size to maintain stability. On the other hand, fully implicit methods can become prohibitively expensive for these types of problems regarding computation time per step and memory requirements. As previously stated, the cost of implicit methods scales with $O({p}^{d})$ for standard FR and $O({p}^{d-1})$ for hybridized FR, where $p$ is the polynomial degree, and $d$ is the dimension of the problem. Hence, the feasibility of employing implicit time-stepping for large-scale computations is limited at high orders.

A more efficient approach involves both explicit and implicit IMEX time-stepping methods. These schemes are able to leverage the stability of implicit schemes for stiff terms while mitigating their cost by using explicit methods for nonstiff terms. To demonstrate their application, consider an ordinary differential equation of the form
\begin{equation}
  u' = f(u) + g(u),
\end{equation}
where $f(u)$ is the nonstiff part of the problem and $g(u)$ is the stiff portion. Consider an implicit $s$-stage diagonally-implicit Runge-Kutta (DIRK) scheme for the stiff region associated with a matrix and vector of coefficients of a Butcher tableau~\cite{butch} given by $\bm A\in\mathbb{R}^{s\times s},~\bm b\in\mathbb{R}^s,~\bm c\in\mathbb{R}^s$. For the nonstiff part, consider an explicit $\sigma=s+1$-stage RK method with respective coefficients given by $\bar{\bm{A}}\in\mathbb{R}^{\sigma\times\sigma},~\bar{\bm{b}}\in\mathbb{R}^\sigma,~\bar{\bm{c}}\in\mathbb{R}^\sigma$. To compensate for the difference in the size of the matrices of coefficients, a first row and a first column of zeros are padded into the implicit tableau. In order to be paired, the implicit and explicit schemes must satisfy $\bar{\bm c} = \left[0~ \bm c\right]^T$. The resulting form of a general IMEX Butcher tableau can be seen in Table~\ref{table:butcherimex}.
\begin{table}[hbtp]
  \centering
  \begin{subtable}{.49\linewidth}\centering{
\begin{tabular}{c|ccccc}
    0     & 0        &  0        & 0 & $\ldots$  & 0 \\
    $c_1$ & 0        & $a_{11}$  & 0 & $\ldots$  & 0 \\ 
    $c_2$ & 0        & $a_{21}$  & $a_{22}$ & $\ldots$  & 0 \\ 
    $\vdots$ & $\vdots$ & $\vdots$ & $\vdots$ & $\ddots$ & $\vdots$\\
    $c_s$    & 0 & $a_{s1}$ & $a_{s2}$  & $\ldots$ & $a_{ss}$ \\ \hline
    & 0 & $b_1$ & $b_2$ & $\ldots$ & $b_s$
\end{tabular}\quad\quad}
\caption{Implicit part}
\label{table:implicitimex}
\end{subtable}
\begin{subtable}{.49\linewidth}\centering{
\begin{tabular}{c|ccccc}
  $\bar c_1$     & 0        &  0        & 0 & $\ldots$  & 0 \\
  $\bar c_2$ & $\bar a_{21}$ & 0  & 0 & $\ldots$  & 0 \\ 
  $\bar c_3$ & $\bar a_{31}$ & $\bar a_{32}$  & 0 & $\ldots$  & 0 \\ 
  $\vdots$ & $\vdots$ & $\vdots$ & $\vdots$ & $\ddots$ & $\vdots$\\
  $\bar c_\sigma$ & $\bar a_{\sigma 1}$ & $\bar a_{\sigma 2}$ & $\bar a_{\sigma 3}$ & $\ldots$ & 0 \\ \hline
  & $\bar b_1$ & $\bar b_2$ & $\bar b_3$ & $\ldots$ & $\bar b_\sigma$
\end{tabular}}
\caption{Explicit part}
\label{table:explicitimex}
\end{subtable}
\caption{General form of Butcher tableaus for IMEX schemes} \label{table:butcherimex}
\end{table}

To advance the solution from time level $n$ to $n+1$ by a time-step $\Delta t$, the first stage is always explicit. Then implicit and explicit solves are alternated, as shown in Algorithm~\ref{alg:imex}~\cite{perssonHighorderSimulationsUsing2011,vermeire2015adaptive}. Using IMEX schemes can yield significantly smaller implicit systems to solve, as it is dedicated to only a fraction of the problem when considering geometry-induced stiffness. In addition, these schemes are linearly stable, maintain the expected orders of accuracy~\cite{vermeire2015adaptive}, and can have superior performance compared to purely explicit and purely implicit methods for LES simulations~\cite{perssonHighorderSimulationsUsing2011}. Furthermore, optimized IMEX methods~\cite{vermeire2021accelerated} can be obtained, resulting in additional speedups. The IMEX approach can be further leveraged by introducing hybridization to solve the implicit portion at each stage. Hence, this section develops an efficient IMEX formulation for geometry-induced stiffness. 

\begin{algorithm}
  \SetAlgoLined
  \DontPrintSemicolon
 Set\
 $\bar R_1 = f(u_{n})$.\;
\For{$i\gets1$ \KwTo $s$}{
  Solve for $R_i$ in $R_i=g(u_i)$, where\;
  \quad $u_i = u_{n} + \Delta t \sum\limits_{j=1}^i a_{i,j} R_j + \Delta t \sum\limits_{j=1}^i \bar{a}_{i+1,j} \bar R_j$.\;
  Evaluate\;
  \quad $\bar R_{i+1} = f(u_i)$.
}
Compute the value at the next time step \;
\quad $u_{n+1} = u_{n} + \Delta t \sum\limits_{j=1}^s b_j R_j + \Delta t \sum\limits_{j=1}^\sigma \bar b_j \bar R_j$.
 \caption{Time integration using an IMEX scheme for one time-step}
 \label{alg:imex}
\end{algorithm}

This section proposes an IMEX formulation by pairing conventional FR and HFR methods to tackle geometry-induced stiffness. While hybridized methods have also been developed in explicit formulations, they require explicit trace definitions~\cite{stanglmeier2016explicit,nehmetallah2020explicit} and nonlinear solvers~\cite{samii2019comparison} and hence the benefits over a standard FR formulation in a general nonlinear problem are still not clear. The explicit form of the FR method is suitable for nonstiff problems. FR methods are locally conservative and have demonstrated potential for modern parallel computer architectures. Hence, we employ it to solve the moderate to large elements associated with lower stiffness in the domain. We introduce hybridization, which is expected to reduce the size of the implicit solver by employing HFR or EFR formulations for the smallest elements associated with the stiff portions of the domain. We refer to the proposed approach as hybridized IMEX methods. These schemes are expected to reduce the cost of a purely implicit method while increasing the allowable time-step size and improving the constrained stability posed by the explicit formulation.

We are interested in integrating the equation
\begin{equation}
  \frac{d \bm u}{dt} = \bm R(\bm u(t)),
  \label{eq:imexode}
\end{equation}
subject to an appropriate initial condition, where $\bm R$ typically contains the divergence of the flux after applying a spatial discretization such as the FR method. In order to integrate this equation for geometry-induced stiffness, an $s$-stage IMEX method with order $q$ to advance a solution from time level $n$ to $n+1$ will be employed according to a modified version of Algorithm~\ref{alg:imex}. 
\begin{figure}[htbp]
  \centering
  \includegraphics[width=0.7\textwidth]{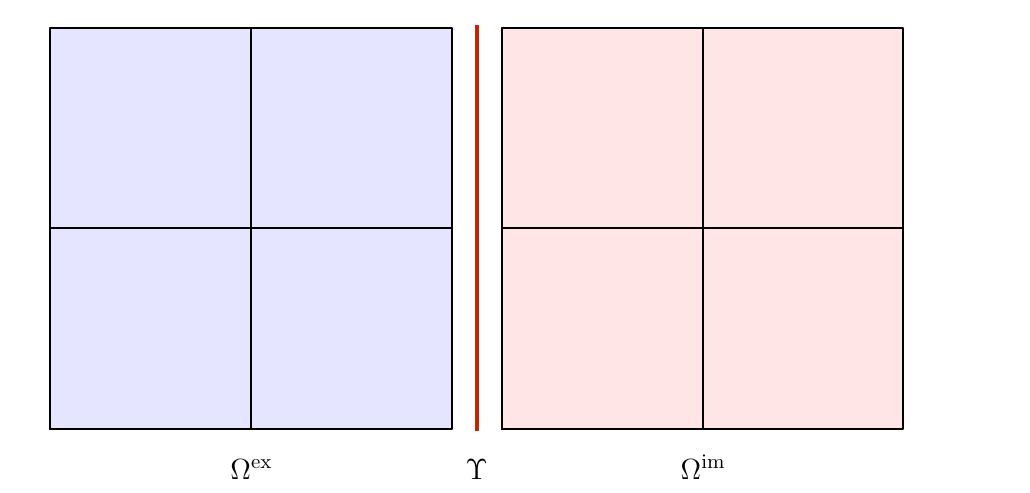}
  \caption{Reference domain partitioning for IMEX schemes}
  \label{fig:imexsketch}
\end{figure}

To this end, consider a computational domain $\Omega$ subdivided into two regions, as shown in Figure~\ref{fig:imexsketch}. The implicit region is denoted $\Omega^{\text{im}}$ and the explicit region $\Omega^{\text{ex}}$. The interface between these two regions is denoted $\Upsilon=\Omega^{\text{ex}}\cap\Omega^{\text{im}}$. The solution and trace vectors can be defined
\begin{equation}
  \bm u = 
  \begin{bmatrix}
  \bm u^{\text{ex}} \\
  \bm u^{\text{im}} 
  \end{bmatrix},
  \quad
  \bm{\hat u} = 
  \bm {\hat u^{\text{im}}},
\end{equation}
respectively, where $\bm u^{\text{im}}\in\mathbb{R}^{N^{\text{im}} N_s},~\bm u^{\text{ex}}\in\mathbb{R}^{N^{\text{ex}} N_s}$, and $\hat{\bm u}^{\text{im}}\in\mathbb{R}^{\hat N N_s}$. $N^{\text{im}},~N^{\text{ex}},\hat N$ are the number of implicit elements, explicit elements, and trace points. The interior solution is found in both explicit and implicit subdomains, but the trace is only defined at $\varepsilon^{h,\text{im}}\setminus\Upsilon$ since $\hat{\bm {u}}^{\text{ex}}=\emptyset$. Hence, we will refer to the trace in the implicit side as $\hat{\bm u}$. After hybridization of the implicit portion and applying the above definitions, the ODE in Equation~\eqref{eq:imexode} can be generalized to a system of the form
\begin{align}
  \frac{d \bm u}{dt} =
  \begin{bmatrix} 
  \bm R^{\text{ex}}(\bm u)\\
  \bm R^{\text{im}}(\bm u, \hat{\bm{u}})
  \end{bmatrix}&\quad\text{in }~\mathcal{T}^h,\\
  \bm G(\bm u, \bm{\hat u}) = \bm 0&\quad\text{in }\varepsilon^{h,\text{im}},\end{align}
\begin{figure}[htbp]
  \centering
  \begin{subfigure}[b]{0.7\textwidth}
    \includegraphics[width=\textwidth]{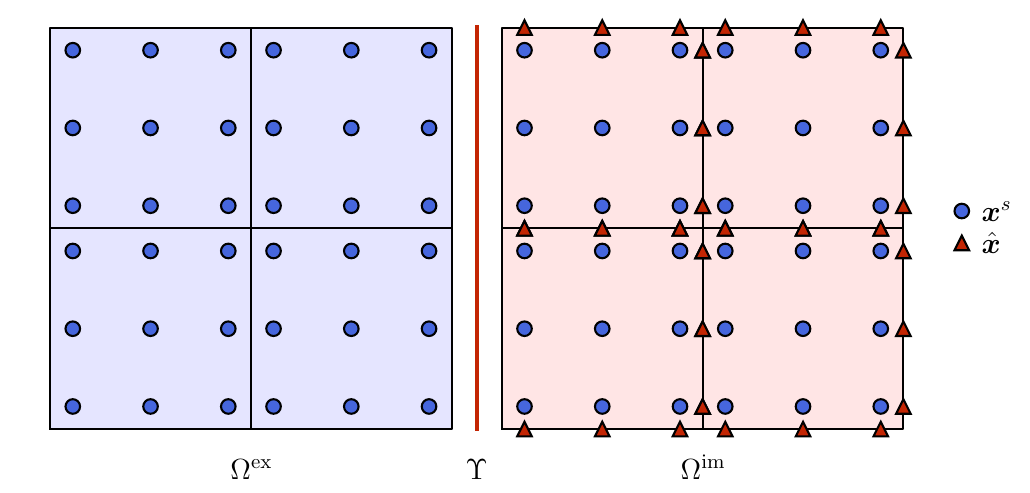}
    \caption{HFR}
  \end{subfigure}
    \begin{subfigure}[b]{0.7\textwidth}
    \includegraphics[width=\textwidth]{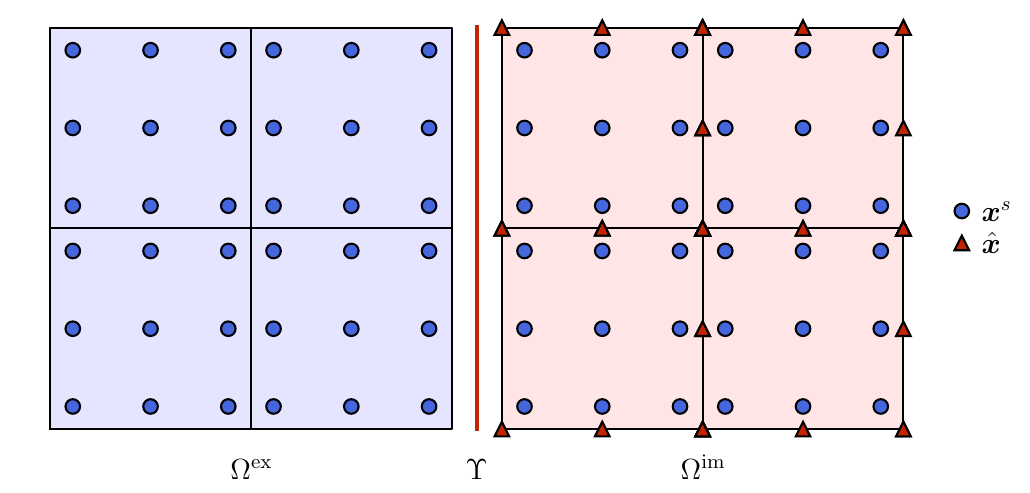}
    \caption{EFR}
    \end{subfigure}
  \caption{Distribution of trace and solution points in this configuration for HFR and EFR IMEX methods for a $p=2$ discretization}
  \label{fig:imexsketchsoltrace}
\end{figure}where $\varepsilon^{h,\text{im}} = \varepsilon^{h} \cap \Omega^{\text{im}}$. The resulting distribution of trace and solution points, $\hat{\bm x}$ and $\bm x^s$, respectively, can be seen in Figure~\ref{fig:imexsketchsoltrace} for HFR and EFR methods. Furthermore, $\bm R^{\text{ex}}(\bm u)$ is associated with the standard FR equations in \eqref{eq:freq}, $\bm R^{\text{im}}(\bm u, \hat{\bm u})$ with the hybridized Equation in~\eqref{eq:hfr1} and $\bm G(\bm u, \hat{\bm u})$ with the transmission conditions in \eqref{eq:hfr2}. Due to the temporally non-overlapping nature of IMEX schemes, the explicit equations can be solved separately from the implicit portion. Applying the IMEX-RK method, the solution at each stage can be found by first computing the value of the explicit solution via 
\begin{equation}
  \bm u_i^{\text{ex}} = \bm u_n^{\text{ex}} + \Delta t \sum_{j=1}^i \bar{a}_{i+1,j} \bm R^{\text{ex}}_j,
  \label{eq:explicitimex}
\end{equation}
where $\bm R^{\text{ex}}_j = \bm R(\bm u_{j-1})$ for $j>1$ and $\bm R^{\text{ex}}_j = \bm R(\bm u_{n})$ for $j=1$. The explicit residual only depends on known values of the solution at stages $i<s$ since $\bar a_{i,j}=0$ for $j\geq i$. Here, the residual is a function of the solution in the explicit subdomain and in the elements in direct contact with the IMEX interface of the implicit portion, which is always known for the required indices. After computing~\eqref{eq:explicitimex}, the implicit solution can be readily obtained via
\begin{subequations}
  \begin{align}
  \bm u_i^{\text{im}} = \bm u_n^{\text{im}} + \Delta t \sum_{j=1}^i a_{i,j} \bm R_j^{\text{im}}(\bm u_j, \hat{\bm u}_j),\\ 
  \bm G(\bm u_i, \bm {\hat u}_i) = 0,\label{eq:galdkjad}
  \end{align}
\end{subequations}
which employs hybridization and Equation~\eqref{eq:galdkjad} represents discrete transmission conditions with the form
\begin{equation}
  \sum_{\bar f\in \varepsilon^{h,\text{im}}_0\setminus \Upsilon} \int_{\bar f} \llbracket \hat{\bm{\mathfrak{F}}}(u_i,\hat u_i)\rrbracket_{\bar f} \phi ds
  + \sum_{\bar f\in \Upsilon} \int_{\bar f} \llbracket\hat{{\bm {\mathfrak F}}}^{\text{FR}}(u_i^{\text{im}},u_i^{\text{ex}})\rrbracket_{\bar f} \phi ds 
  + \sum_{\bar f\in \varepsilon^{h,\text{im}}_\partial}  \int_{\bar f}  {\mathfrak{F}}^{\text{BC}}_{\bar f} \phi ds = 0,\label{eq:hfrimex}
\end{equation}
where the typical transmission equations of Equation~\eqref{eq:hfr2} have been augmented with an interface condition to weakly enforce conservation along the IMEX interface and hence globally in the domain. Then, the solution at the next step can be found by
\begin{align}
  \bm u_{n+1} = \bm u_{n} + \Delta t \sum_{j=1}^s b_j \bm R^{\text{im}}_j + \Delta t \sum_{j=1}^\sigma \bar b_j \bm R^{\text{ex}}_j.
\end{align}
At the IMEX interface, we apply the standard FR fluxes and introduce them as boundary conditions for the hybridized portion. Consequently, at the interface, there is no trace definition. The local conservation property of the FR method~\cite{huynhFluxReconstructionApproach2007} enables the use of this approach to both HFR and EFR methods. Recall that our EFR implementation uses discontinuous traces at the boundaries throughout this work. In the proposed hybridized formulation, discontinuous traces are also used at the IMEX interface for both methods. This allows for using paired HFR-FR and EFR-FR methods to tackle geometry-induced stiffness. The nonstiff portion will retain its local conservation properties, and the stiff portion will be globally conservative for the EFR method and locally conservative for the HFR method. Specifically, for the HFR methods with discontinuous traces, the transmission conditions reduce to a pointwise conservation statement, which can be shown to yield the same definition of the traces in linear advection throughout the domain and recovers a standard FR IMEX method in these cases only. The proof is trivial and is omitted for brevity.

\section{Numerical Examples}
This section presents a series of numerical examples to showcase the benefits of using our proposed hybridized IMEX formulations to tackle geometrical stiffness in two and three dimensions. We will first present verification of our approach via linear advection and then demonstrate capabilities via nonlinear numerical examples. To evaluate geometry-induced stiffness, the following elementwise stiffness indicator is used
\begin{equation}
 \operatorname{E}_s = \frac{|\Omega_k|}{|\partial\Omega_k|},
\end{equation}
where $|\Omega_k|$ is the volume of the element and $|\partial\Omega_k|$ is the sum of the face areas. Hence, the indicator has dimensions of length. For elements with high aspect ratios as well as for very small elements, the sum of the face areas is significantly larger than its volume, which means that elements with high numerical stiffness will have a high value of $\operatorname{E}_s$. The use of a user-provided cutoff $\operatorname{E}_s$ determines the implicit ratio, which we define to be the number of implicit elements $N^{\text{im}}$ in relation to the total number of elements $N$ in the domain
\begin{equation}
\operatorname{IF} = \frac{N^{\text{im}}}{N},
\end{equation}
so that a higher cutoff value of $\operatorname{E}_s$ yields more implicit elements. For simplicity, this work will refer to FR, HFR, and EFR methods as IMEX discretizations with implicit portions solved using the FR, HFR, and EFR methods, respectively. For all runs, the explicit portion uses a standard FR discretization. In the following problems, IMEX$_{s,q}$ methods with $s$ stages and $q$ order are considered. 
\subsection{Verification}

\begin{figure}[htbp]
  \centering
  \begin{subfigure}[b]{0.49\textwidth}
    \includegraphics[width=\textwidth]{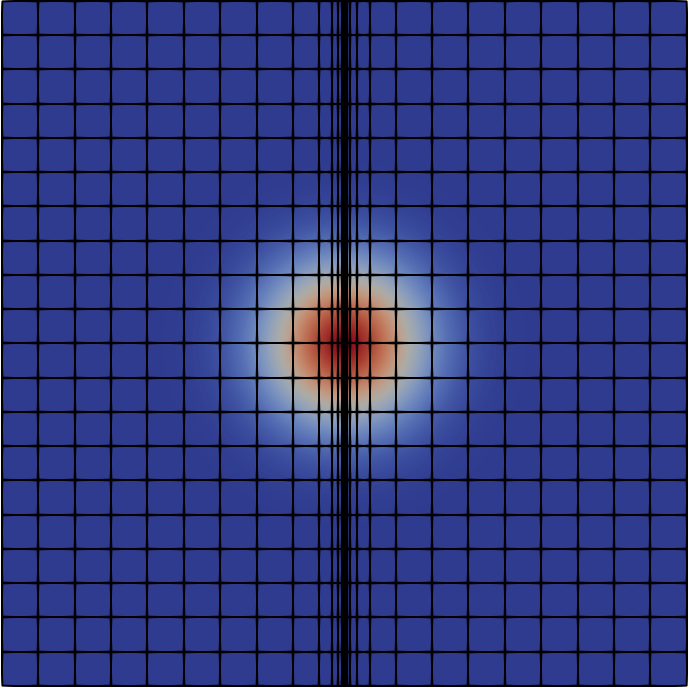}
    \caption{Final solution}
  \end{subfigure}
    \begin{subfigure}[b]{0.49\textwidth}
    \includegraphics[width=\textwidth]{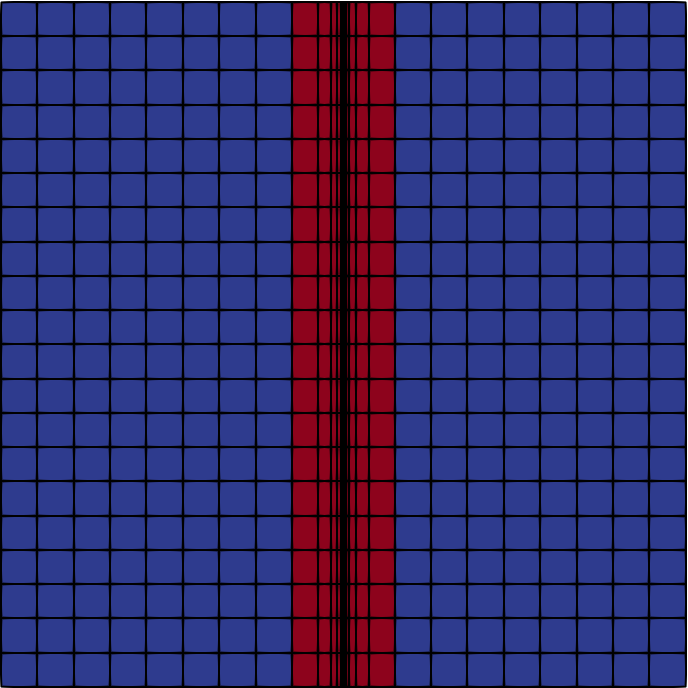}
    \caption{Implicit-explicit distribution}
    \end{subfigure}
  \caption{Distribution of explicit and implicit elements for the IMEX verification of linear advection}
  \label{fig:advectioncontours}
\end{figure}

We perform verification using linear advection by considering a periodic square domain. The domain is split into quadrilateral elements. Along the center, a band of seven layers of stretched elements is placed with a stretching ratio equal to 2, as shown in Figure~\ref{fig:advectioncontours}. Within this band, elements are flagged as implicit, whereas uniform elements away from this section remain explicit. The implicit portion of the domain is solved using a hybridized form, and the explicit portion uses a standard FR method. Both implicit and explicit regions use a solution polynomial degree of seven to reduce contamination arising from spatial error. The initial condition is a Gaussian profile
\begin{equation}
  u(\bm x, 0) = \exp\left(-\frac{1}{20}[(x-x_c)^2 + (y-y_c)^2]\right),
\end{equation}  
where $(x_c,~y_c)=(10,~10)$ is the center coordinate of the domain of size $20\times20$. After one convective time $t_c=20$, the $L_2$-norm of the error is computed. Results of the $L_2$ norm of the error against the exact solution are shown in Table~\ref{tab:ooa-adv-imex} for the IMEX$_{3,2}$~\cite{ascherImplicitexplicitRungeKuttaMethods1997}, IMEX$_{5,3}$~\cite{kanevskyApplicationImplicitExplicit2007} and AIMEX$_{10,2}$ method~\cite{vermeire2021accelerated} with Butcher tableaus included in the appendix for the IMEX methods and are available for AIMEX as supplementary material in~\cite{vermeire2021accelerated}. Due to the temporal error dominating the $L_2$ norm, results for the EFR and HFR methods differ only beyond single precision, and no difference can be observed in the tabulated values. In addition, results for the approach with HFR are equivalent to solving the problem using FR everywhere in the domain, as expected for linear advection. Hence, only the HFR approach is shown. The second and third orders of temporal accuracy are recovered for the three considered IMEX schemes. Of the three considered IMEX methods, AIMEX$_{10,2}$ provides relatively lower error than IMEX$_{3,2}$ for a given time-step size $\Delta t$, allowing for a larger time step due to the optimized explicit stability polynomial with only two implicit solves. While the IMEX$_{5,3}$ method is more accurate than the other two, it requires five implicit solves per time step. Thus, we use the optimized AIMEX$_{10,2}$ method for the rest of this work.

\begin{table}[htbp]
\centering
\caption{Convergence table showing the $L_2$-norm of the solution error and the order of accuracy for linear advection using multiple IMEX schemes}
\label{tab:ooa-adv-imex}
\begin{tabular}{cccccccc}
\hline
Scheme                          & $\Delta t$ &  & HFR      & Order  &  & EFR      & Order  \\ \hline
\multirow{4}{*}{IMEX$_{3,2}$}   & 0.02       &  & 6.78$\times10^{-6}$ & -    &  & 6.78$\times10^{-6}$ & -    \\
                                & 0.01       &  & 1.67$\times10^{-6}$ & 2.02 &  & 1.67$\times10^{-6}$ & 2.02 \\
                                & 0.005      &  & 4.17$\times10^{-7}$ & 2.00 &  & 4.17$\times10^{-7}$ & 2.00 \\
                                & 0.0025     &  & 1.04$\times10^{-7}$ & 2.00 &  & 1.04$\times10^{-7}$ & 2.00 \\ \hline
\multirow{4}{*}{IMEX$_{5,3}$}   & 0.02       &  & 2.19$\times10^{-6}$ & -    &  & 2.19$\times10^{-6}$ & -    \\
                                & 0.01       &  & 2.73$\times10^{-7}$ & 3.00 &  & 2.73$\times10^{-7}$ & 3.00 \\
                                & 0.005      &  & 3.41$\times10^{-8}$ & 3.00 &  & 3.41$\times10^{-8}$ & 3.00 \\
                                & 0.0025     &  & 4.28$\times10^{-9}$ & 3.00 &  & 4.28$\times10^{-9}$ & 3.00 \\ \hline
\multirow{4}{*}{AIMEX$_{10,2}$} & 0.02       &  & 4.28$\times10^{-6}$ & -    &  & 1.83$\times10^{-6}$ & -    \\
                                & 0.01       &  & 4.60$\times10^{-7}$ & 1.99 &  & 4.60$\times10^{-7}$ & 1.99 \\
                                & 0.005      &  & 1.15$\times10^{-7}$ & 2.00 &  & 1.15$\times10^{-7}$ & 2.00 \\
                                & 0.0025     &  & 2.89$\times10^{-8}$ & 2.00 &  & 2.89$\times10^{-8}$ & 2.00 \\ \hline 
\end{tabular}
\end{table}

\subsection{Laminar Flow over a Circular Cylinder}
\begin{figure}[htbp]
  \centering
   \includegraphics[width=\textwidth]{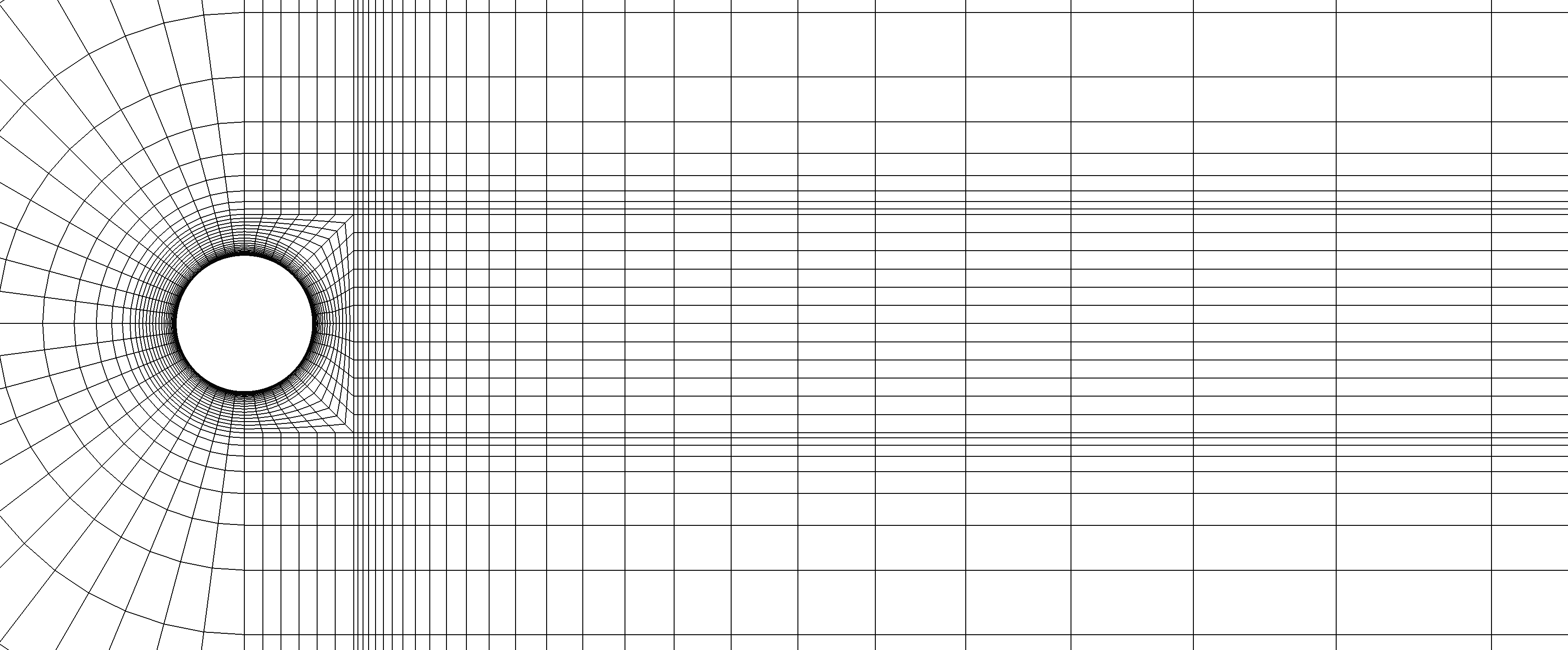}
  \caption{Mesh for the laminar cylinder case at $\operatorname{Re}=150$ consisting of 3090 quadrilateral elements}
  \label{fig:cyl150mesh}
\end{figure}

\begin{figure}[htbp]
  \centering
   \includegraphics[width=0.65\textwidth]{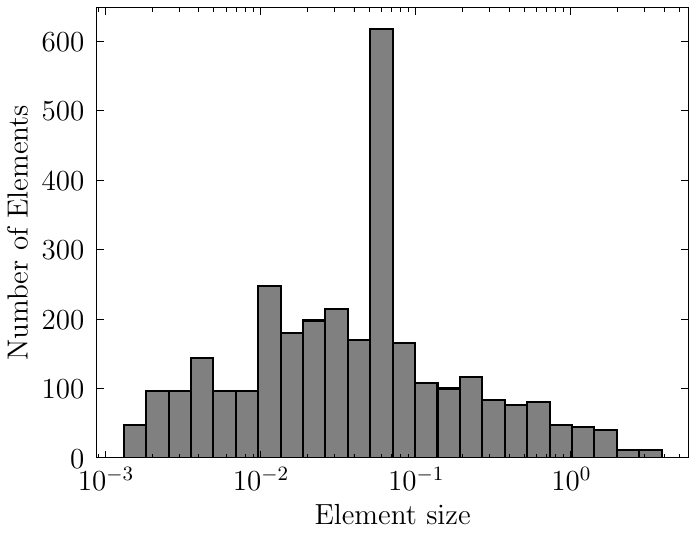}
  \caption{Distribution of element sizes for the laminar cylinder case}
  \label{fig:cylhistogram}
\end{figure}

In this section, we study flow over a cylinder at $\operatorname{Re}=150$, based on the cylinder diameter $D$, to simulate unsteady vortex shedding. This problem belongs to the laminar regime. Hence, a two-dimensional approach is suitable for this simulation. Here, we want to validate our implementation of hybridized IMEX discretizations. We make use of a computational domain divided into 3090 quadrilateral elements refined toward the cylinder walls to capture the gradients due to the boundary layer, as shown in Figure~\ref{fig:cyl150mesh}. We use this case as a baseline problem to analyze the performance of IMEX Hybridized-FR schemes. We choose to run this problem at Mach number $\operatorname{M}=0.1$ to compare against the reference data of Cagnone~\cite{cagnoneStableInterfaceElement2012}. In this case, the resolution within the boundary layer is increased to test the benefit of these methods in stiff regions. In Figure~\ref{fig:cylhistogram}, a histogram containing the distribution of element sizes is provided. We observe disparate element sizes with differences of up to three orders of magnitude between the largest and the smallest elements. 

\begin{figure}[htbp]
  \centering
  \begin{subfigure}[b]{0.49\textwidth}
    \includegraphics[width=\textwidth]{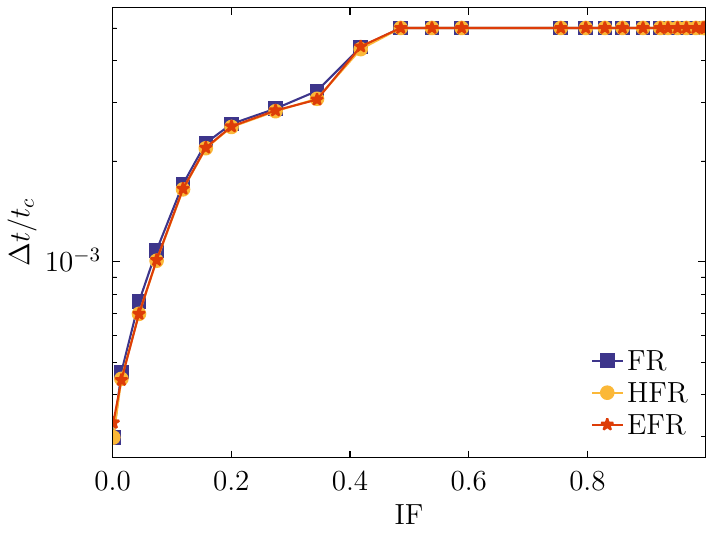}
    \caption{$p=1$}
  \end{subfigure}
    \begin{subfigure}[b]{0.49\textwidth}
    \includegraphics[width=\textwidth]{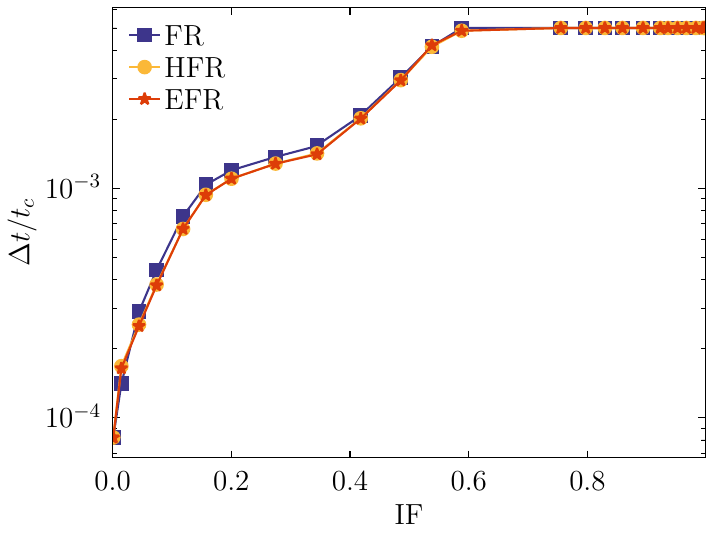}
    \caption{$p=2$}
  \end{subfigure}
  \begin{subfigure}[b]{0.49\textwidth}
    \includegraphics[width=\textwidth]{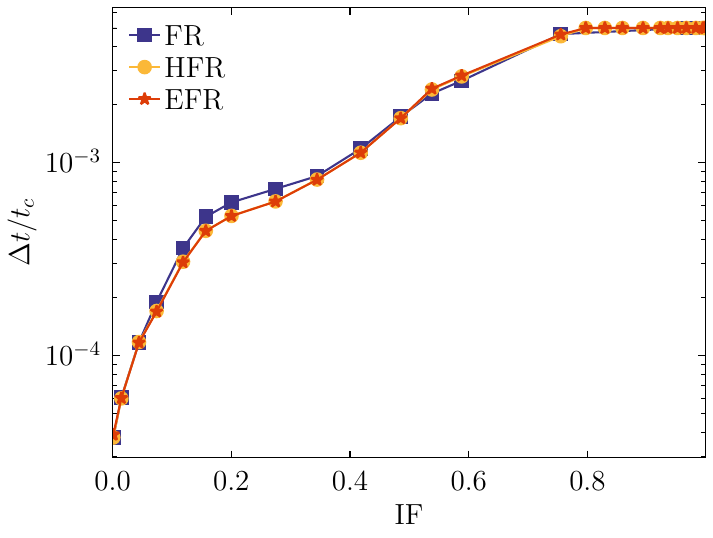}
    \caption{$p=3$}
  \end{subfigure}
    \begin{subfigure}[b]{0.49\textwidth}
    \includegraphics[width=\textwidth]{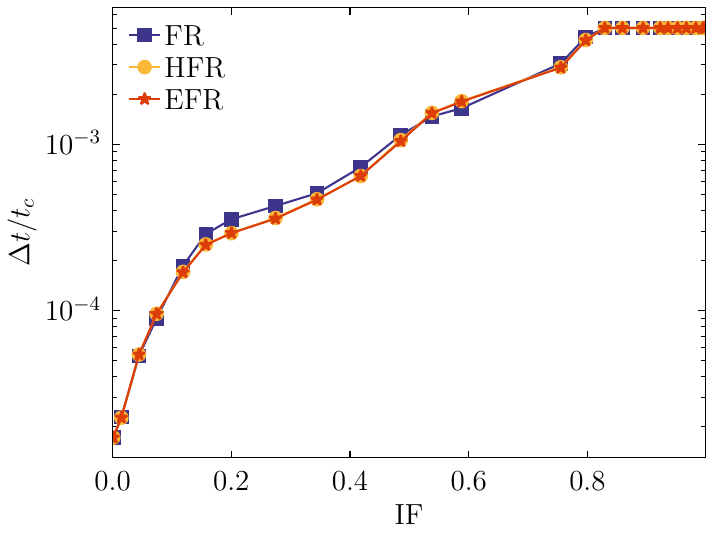}
    \caption{$p=4$}
  \end{subfigure}
  \caption{Maximum stable time-step size for multiple implicit fractions for the AIMEX$_{10,2}$ scheme}
  \label{fig:dtmaxaimex2}
\end{figure}

\begin{figure}[htbp]
  \centering
  \begin{subfigure}[b]{0.49\textwidth}
    \includegraphics[width=\textwidth]{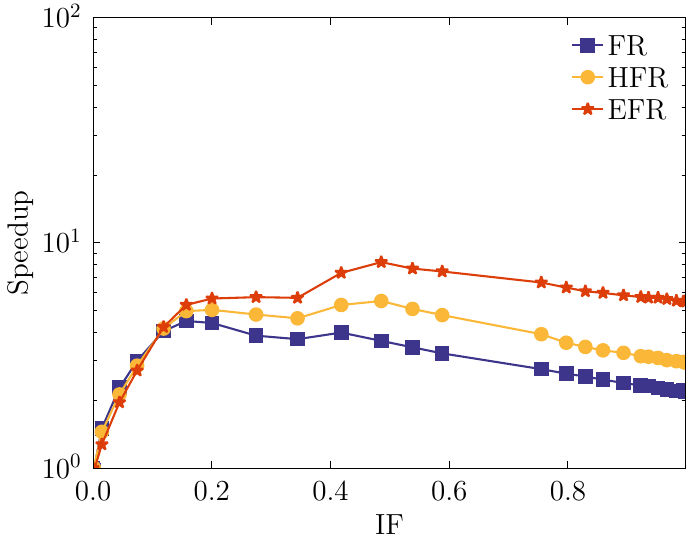}
    \caption{$p=1$}
  \end{subfigure}
    \begin{subfigure}[b]{0.49\textwidth}
    \includegraphics[width=\textwidth]{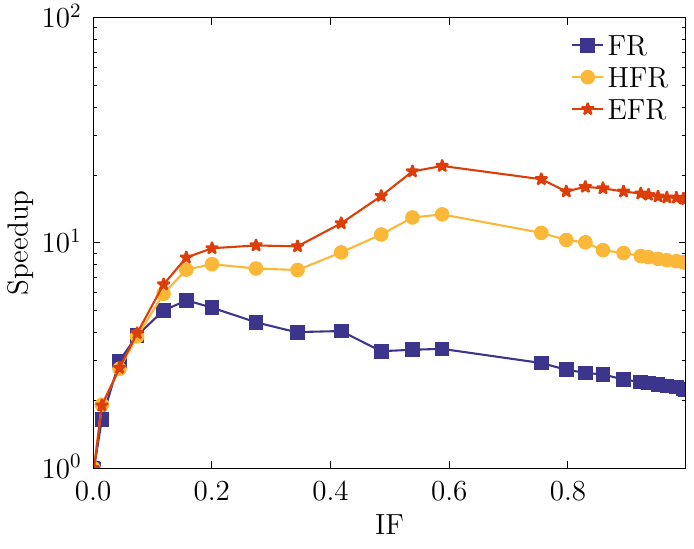}
    \caption{$p=2$}
  \end{subfigure}
  \begin{subfigure}[b]{0.49\textwidth}
    \includegraphics[width=\textwidth]{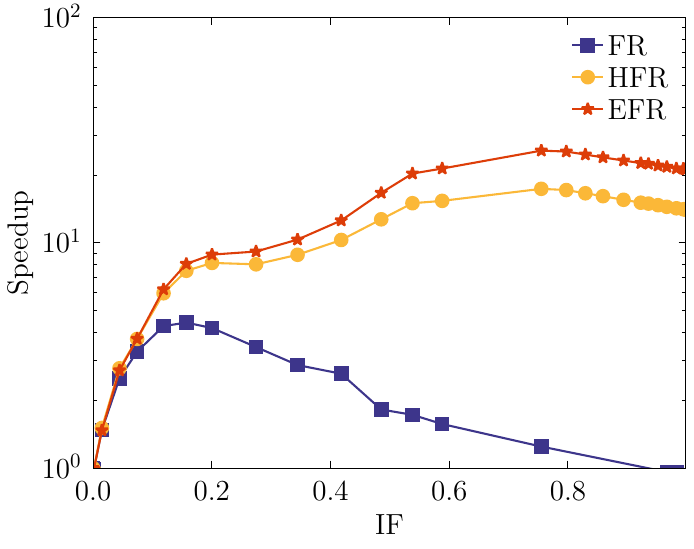}
    \caption{$p=3$}
  \end{subfigure}
    \begin{subfigure}[b]{0.49\textwidth}
    \includegraphics[width=\textwidth]{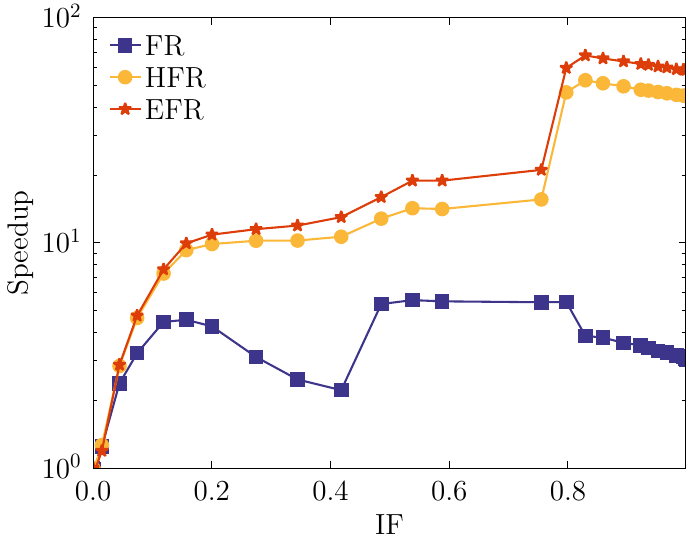}
    \caption{$p=4$}
  \end{subfigure}
  \caption{Speedup factors for multiple implicit fractions for the AIMEX$_{10,2}$ scheme}
  \label{fig:speedupsaimex2}
\end{figure}

By increasing the cutoff element size, the maximum allowable time-step size is expected to increase. The maximum stable $\Delta t$ are computed for a range of implicit factors using a simple bisection algorithm and shown in Figure~\ref{fig:dtmaxaimex2}. We capped the maximum time-step size to $t_c/200$ to maintain accuracy. This value is typically achieved at $\operatorname{IF}>0.5$, a relatively large implicit factor that requires significant computational resources in industrial-scale problems, as we will discuss in the last numerical example. Our main interest is focused on regions with moderately low implicit factors. Note from this plot that the value of $\Delta t_{\max}$ is similar for the FR, HFR, and EFR methods, with a few differences due to the shape of the stability polynomial in AIMEX methods optimized for the FR method only. We compute speedups against FR explicit formulations with $\operatorname{IF}=0$. Speedups near two orders of magnitudes for $p=4$ simulations are observed at higher IF values in Figure~\ref{fig:speedupsaimex2}. These speedups were computed using serial simulations. We observed hybridized methods to yield up to 2.5 times faster results than an FR method at $\operatorname{IF}=0.2$. This translates to 10 times faster results compared to an explicit approach while utilizing significantly less memory than a fully implicit method. Furthermore, an implicit factor $\operatorname{IF}=0.2$ is chosen, which is expected to provide optimal speedup factors without using significant memory.


\begin{figure}[htbp]
  \centering
  \begin{subfigure}[b]{0.49\textwidth}
    \includegraphics[width=\textwidth]{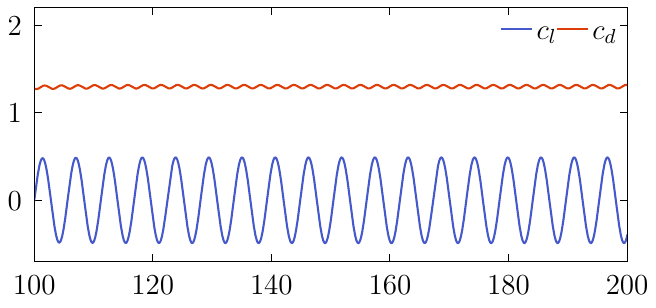}
    \caption{FR, $p=1$}
  \end{subfigure}
    \begin{subfigure}[b]{0.49\textwidth}
    \includegraphics[width=\textwidth]{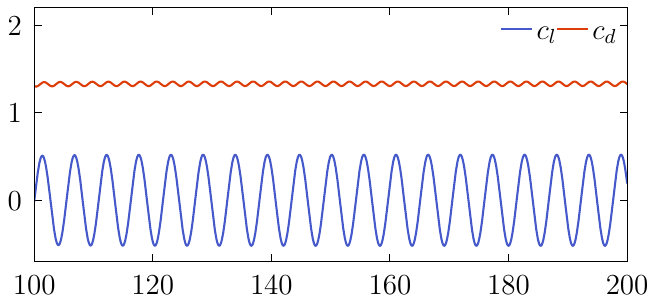}
    \caption{FR, $p=4$}
  \end{subfigure}
  \begin{subfigure}[b]{0.49\textwidth}
    \includegraphics[width=\textwidth]{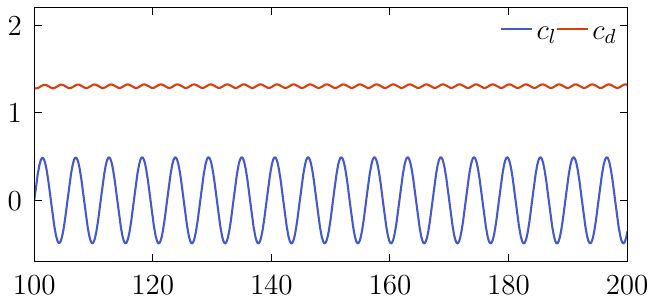}
    \caption{HFR, $p=1$}
  \end{subfigure}
    \begin{subfigure}[b]{0.49\textwidth}
    \includegraphics[width=\textwidth]{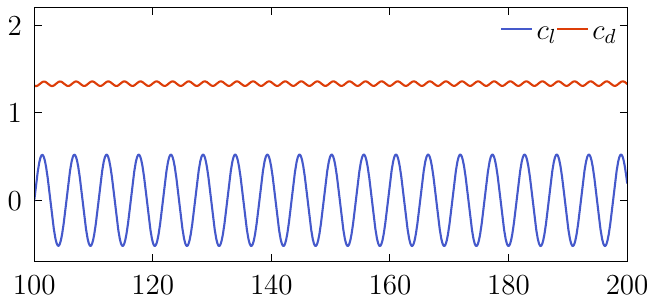}
    \caption{HFR, $p=4$}
  \end{subfigure}
  \begin{subfigure}[b]{0.49\textwidth}
    \includegraphics[width=\textwidth]{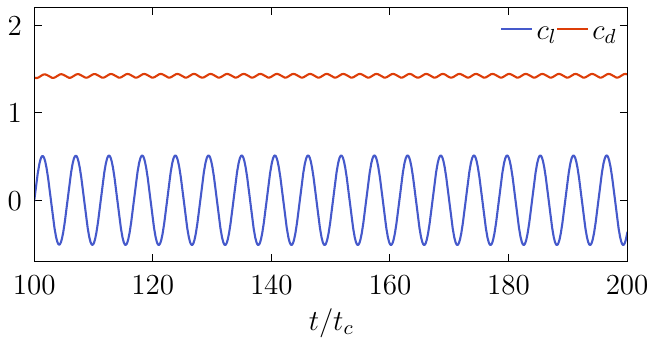}
    \caption{EFR, $p=1$}
  \end{subfigure}
    \begin{subfigure}[b]{0.49\textwidth}
    \includegraphics[width=\textwidth]{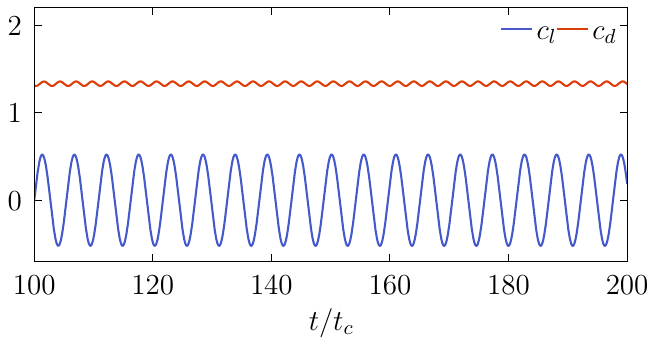}
    \caption{EFR, $p=4$}
  \end{subfigure}
  \caption{Evolution of drag and lift coefficients for the cylinder at $\operatorname{Re}=150$}
  \label{fig:cylinder150clcdimex}
\end{figure}

We perform a series of simulations with polynomial degrees $p=1$ to $p=4$ for 200 convective times $t_c=U_\infty/D$. The evolution of the lift and drag coefficients for the $p=1$ and $p=4$ HFR and EFR schemes is shown in Figure~\ref{fig:cylinder150clcdimex}. We observe a periodic sinusoidal wave associated with the characteristic vortex shedding of this problem. By looking at this plot, it is clear that the frequency of the low-order simulation is different from the $p=4$ results. However, no significant difference is observed between the hybridized and standard IMEX results. More specifically, the Strouhal number for this problem converges to $0.1844$ after $p=3$. $p=1$ is heavily underresolved and underpredicts this result, as shown in Table~\ref{tab:cylinder150results}. The converged results represent less than 0.5\% relative error compared to the numerical results of Cagnone~\cite{cagnoneStableInterfaceElement2012} and less than 1\% compared to the experimental data. Results from the hybridized methods agree with the FR results, with differences of less than 0.2\% in all quantities. Hence, the proposed hybridized IMEX formulations behave similarly to FR, especially at higher orders. The performance benefit of this approach is significant against explicit FR methods. We now evaluate them in a three-dimensional cylinder case in the next section.

\begin{table}[htbp]
\centering
\caption{Summary of results for the cylinder at $\operatorname{Re}=150$ for dual and single scheme IMEX methods}
\label{tab:cylinder150results}
\begin{tabular}{llllll}
\hline
Implicit Scheme      & $p$                & $\bar c_d$ & $\Delta c_d$ & $\Delta c_l$ & $\operatorname{St}$ \\ \hline
\multirow{4}{*}{FR}  & 1                  & 1.2962     & 0.0205       & 0.4879       & 0.1785              \\
                     & 2                  & 1.3271     & 0.0257       & 0.5178       & 0.1842              \\
                     & 3                  & 1.3289     & 0.0258       & 0.5189       & 0.1844              \\
                     & 4                  & 1.3295     & 0.0258       & 0.5192       & 0.1844              \\ \hline
\multirow{4}{*}{HFR} & 1                  & 1.3020     & 0.0206       & 0.4901       & 0.1787              \\
                     & 2                  & 1.3277     & 0.0258       & 0.5187       & 0.1842              \\
                     & 3                  & 1.3267     & 0.0256       & 0.5161       & 0.1844              \\
                     & 4                  & 1.3312     & 0.0259       & 0.5211       & 0.1844              \\ \hline
\multirow{4}{*}{EFR} & 1                  & 1.4207     & 0.0212       & 0.5099       & 0.1787              \\
                     & 2                  & 1.3280     & 0.0258       & 0.5188       & 0.1842              \\
                     & 3                  & 1.3285     & 0.0258       & 0.5178       & 0.1844              \\
                     & 4                  & 1.3309     & 0.0259       & 0.5207       & 0.1844              \\ \hline
\multicolumn{2}{l}{Cagnone~\cite{cagnonePadaptiveLCPFormulation2013}, $p=4$}& 1.3246     & 0.0258       & 0.5166       & 0.1836              \\
\multicolumn{2}{l}{Inoue~\cite{inoueSoundGenerationTwodimensional2002}}     & 1.3200     & 0.0260       & 0.5200       & 0.1830              \\ \hline
\end{tabular}
\end{table}

\subsection{Turbulent Flow over a Circular Cylinder}
\begin{figure}[htbp]
  \centering
   \includegraphics[width=\textwidth]{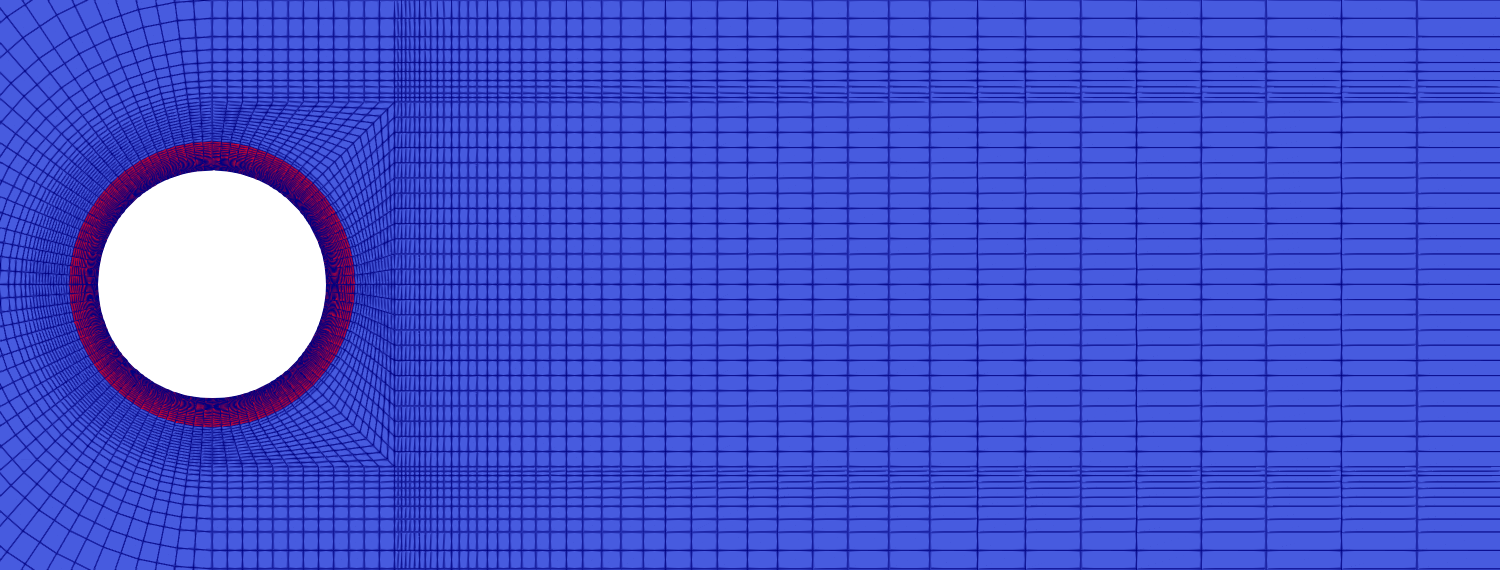}
  \caption{Distribution of implicit (red) and explicit (blue) elements in the computational domain for the turbulent cylinder case}
  \label{fig:cyl1000mesh}
\end{figure}

Three-dimensional flow over a cylinder at $\operatorname{Re}=1000$ is simulated in this section. The computational grid is composed of $37080$ hexahedral elements. Along the spanwise direction, a length of $L_z=2\pi$ is used, a few units over the minimum length required to resolve the three-dimensionality, which is $4D$~\cite{lei2001spanwise}. This length is divided using a grid spacing of $\Delta z = \pi/6D$ in the streamwise direction, which results in 40 layers of elements. Beyond the cylinder, the boundaries were placed at a distance of $40D$ downstream to minimize the effects of boundary conditions. In the previous section, we observed that an implicit factor close to 0.2 performs well in the 2D case and will also be used in this problem. In reality, large-scale problems demand an overwhelming amount of memory, so high implicit factors or fully implicit methods require availability of a vast amount of resources at high orders. Similarly, the AIMEX$_{10,2}$ method is employed here for time integration with $\Delta t/t_c=2.6\times10^{-2}$, which is decreased by half per unit increase in polynomial degree. We distribute implicit and explicit elements as shown in Figure~\ref{fig:cyl1000mesh}, and apply hybridized methods on the implicit portions of the domain. We converge our implicit residuals to a tolerance of $10^{-6}$. The smaller elements in the vicinity of the cylinder walls are flagged as implicit, and larger elements away from the walls are flagged as explicit. 

\begin{figure}[htbp]
  \centering
  \begin{subfigure}[b]{0.49\textwidth}
    \includegraphics[width=\textwidth]{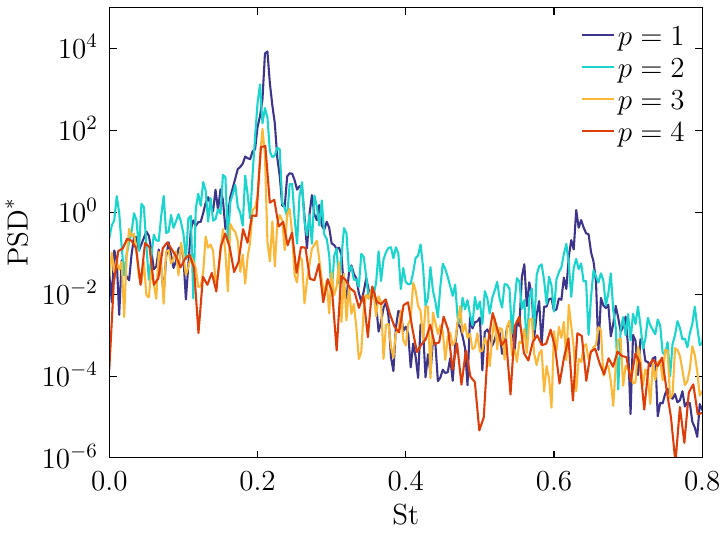}
    \caption{HFR}
  \end{subfigure}
    \begin{subfigure}[b]{0.49\textwidth}
    \includegraphics[width=\textwidth]{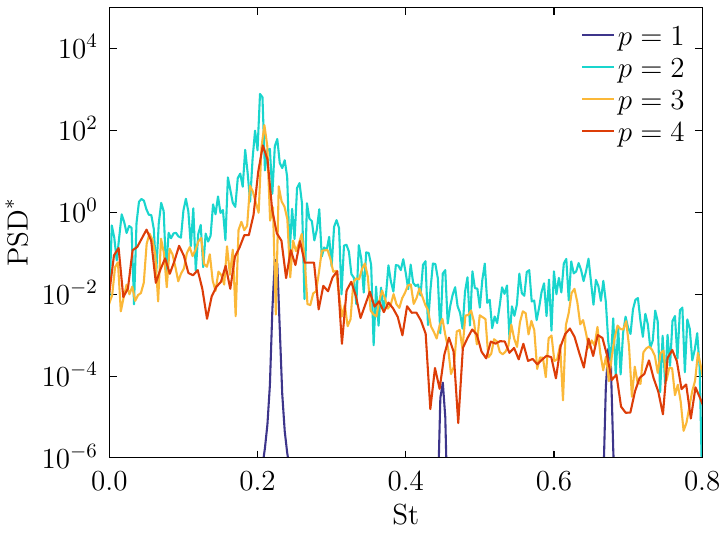}
    \caption{EFR}
  \end{subfigure}
  \caption{Spectra of the unsteady lift coefficient fluctuations for the turbulent cylinder problem}
  \label{fig:clspectraimex}
\end{figure}

We ran this simulation for 200 convective times and averaged the statistics for the last 100$t_c$ to remove the initial transient effects. Results are shown in Table~\ref{tab:cyl1000resultsimex}, where the drag coefficient $C_D$, lift coefficient fluctuations $C_{L,rms}$, and the Strouhal number $\operatorname{St}$ are shown. We compare our results against the data of~\cite{zhaoDirectNumericalSimulation2009}, which provided reference values for a coarse and a fine problem. Relative convergence of the $C_D$ can be observed as the order is increased, which is within 1\% of the coarse results and 15\% of the reference data. The mean lift fluctuations are closer for the HFR method than for the EFR method, which is known to introduce additional error. Overall, results converge to the reference data. The Strouhal number was captured well for the EFR method at $p=4$. In the case of $p=1$, the EFR method did not transition, which caused a significant discrepancy with the reference St as opposed to the other values. The spectra of the $C_L$ signal are shown in~\ref{fig:clspectraimex}, where the convergence to 0.21 is seen as the order is increased. 
\begin{table}[htbp]
\centering
\caption{Summary of results for the turbulent cylinder case}
\label{tab:cyl1000resultsimex}
\begin{tabular}{lrrrr}
\hline
\multicolumn{1}{c}{Scheme} & $p$ & $\bar C_D$ & $C_{L,rms}$ & $\operatorname{St}$ \\ \hline
\multirow{4}{*}{HFR}       & 1   & 2.152     & 1.046      & 0.210               \\
                           & 2   & 1.156     & 0.640      & 0.204               \\
                           & 3   & 1.009     & 0.361      & 0.207               \\
                           & 4   & 0.998     & 0.333      & 0.209               \\ \hline
\multirow{4}{*}{EFR}       & 1   & 2.907     & 0.052      & 0.225               \\
                           & 2   & 1.151     & 0.656      & 0.205               \\
                           & 3   & 1.012     & 0.379      & 0.209               \\
                           & 4   & 1.000     & 0.340      & 0.210               \\ \hline
\multicolumn{2}{l}{Zhao et al. (coarse)~\cite{zhaoDirectNumericalSimulation2009}} & 1.092     & 0.310      & 0.210               \\
\multicolumn{2}{l}{Zhao et al. (fine)~\cite{zhaoDirectNumericalSimulation2009}}   & 1.170      & 0.335      & 0.210                \\ \hline
\end{tabular}
\end{table}

\begin{figure}[htbp]
  \centering
    \begin{subfigure}[b]{0.75\textwidth}
    \includegraphics[width=\textwidth,trim={1cm 1cm 2cm 1cm},clip]{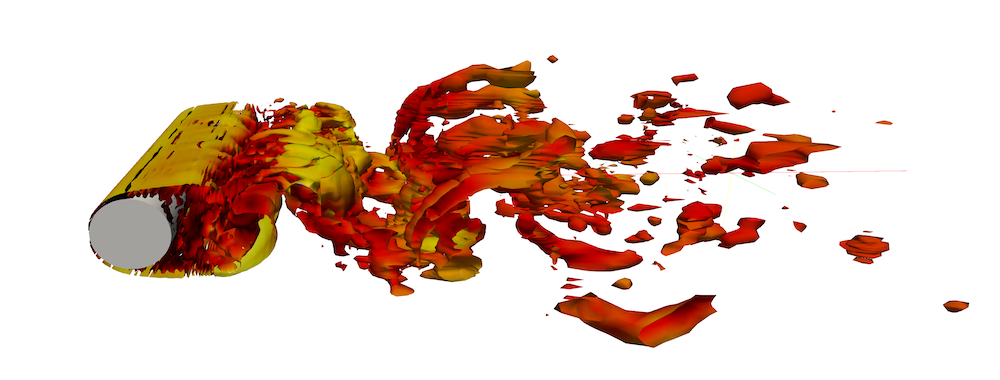}
    \caption{HFR, $p=2$}
  \end{subfigure}
  \begin{subfigure}[b]{0.75\textwidth}
    \includegraphics[width=\textwidth,trim={1cm 1cm 2cm 1cm},clip]{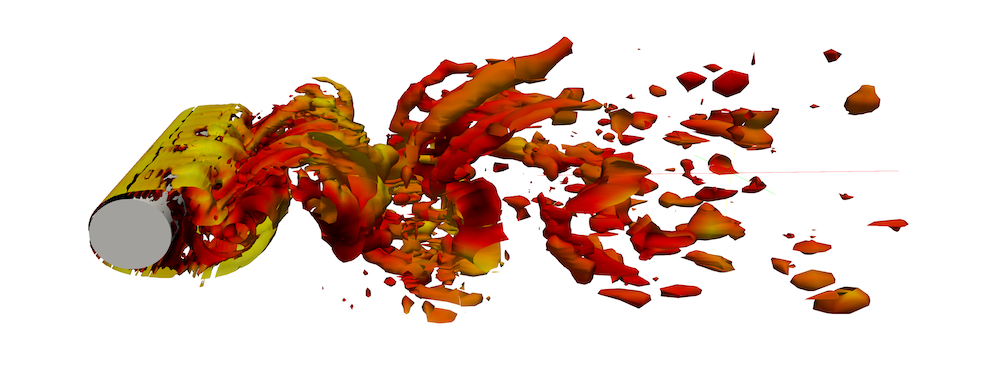}
    \caption{EFR, $p=2$}
  \end{subfigure}
      \begin{subfigure}[b]{0.75\textwidth}
    \includegraphics[width=\textwidth,trim={1cm 1cm 2cm 1cm},clip]{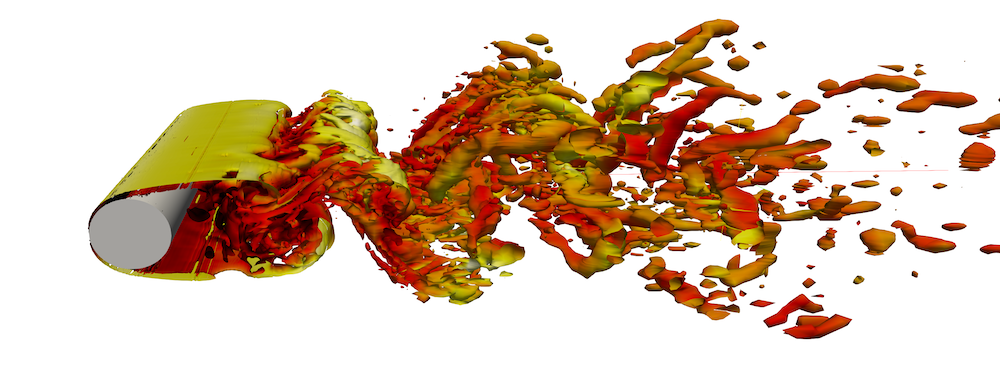}
    \caption{HFR, $p=4$}
  \end{subfigure}
  \begin{subfigure}[b]{0.75\textwidth}
    \includegraphics[width=\textwidth,trim={1cm 1cm 2cm 1cm},clip]{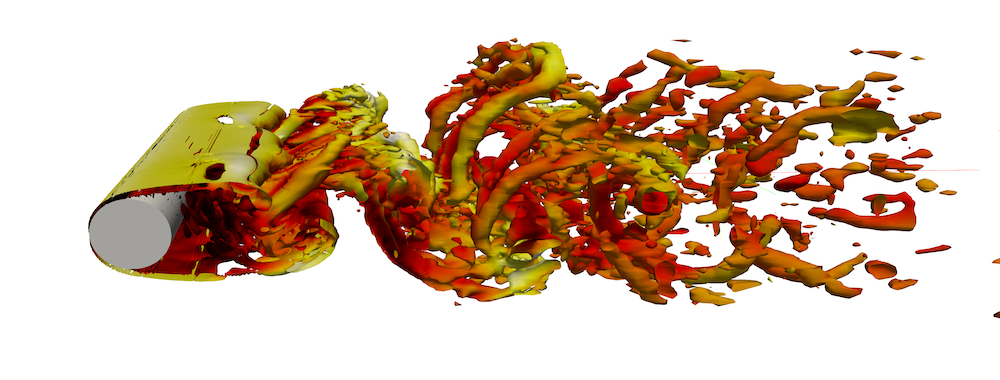}
    \caption{EFR, $p=4$}
  \end{subfigure}
  \caption{Side view of Q-criterion contours for the turbulent cylinder problem}
  \label{fig:qcriterionturbcylside}
\end{figure}

Performance results for this problem demonstrate that hybridized methods are also suitable for three-dimensional problems. Results are shown for the time spent on the global solves $t_G$, local solves $t_L$, block Jacobian computations and implicit matrix assembly $t_J$, and right-hand-side residual computations in the implicit and explicit portions ($t_{R^{\text{im}}}$ and $t_{R^{\text{ex}}}$, respectively), which together add to the overall wall-clock time $t_w$. We show speedups against standard FR IMEX schemes $t_w/t_{w}^{\text{FR}}$ and against explicit runs at $\operatorname{IF}=0$, $t_w/t_{w}^{\text{ex}}$ with the AIMEX tableaus. Results are tabulated in Table~\ref{tab:cyl1000performanceimex}. The timing results are computed for 100 time steps. The time to assemble the Jacobian matrix is associated with a single call, as it was updated every 100 time steps to reduce overhead. These results are computed on 2.4GHz AMD Rome 7532 CPUs using 64 cores. The time spent on the global solutions takes a significantly higher proportion of the total FR computations due to the large size of these systems. For hybridized methods, the solution of the local problems adds overhead, which has linear scalability and can be done efficiently. We observed speedup factors around 2.5 against IMEX and explicit FR using the EFR method, which consistently provided faster runs than FR for $p>1$. The HFR method, however, is only a benefit at $p=4$. This is consistent with the increased number of trace unknowns appearing in hexahedral elements. 
\begin{table}[htbp]
\centering
\caption{Summary of performance metrics for the turbulent cylinder case for 100 time steps}
\label{tab:cyl1000performanceimex}
\begin{tabular}{lrrrrrrrrr}
\hline
Scheme & $p$ & $t_{G}$ & $t_{L}$               & $t_{J}$ & $t_{R^{\text{im}}}$ & $t_{R^{\text{ex}}}$ & $t_w$   & $t_w/t_{w}^{\text{FR}}$ & $t_w/t_{w}^{\text{ex}}$ \\ \hline
FR     & 1   & 187.18  & \multicolumn{1}{c}{-} & 3.10    & 3.72         & 26.39        & 217.33  & 1.00           & 0.76                  \\
       & 2   & 232.34  & \multicolumn{1}{c}{-} & 16.04   & 4.46         & 51.55        & 288.44  & 1.00           & 1.43                  \\
       & 3   & 599.79  & \multicolumn{1}{c}{-} & 85.16   & 8.65         & 116.29       & 724.92  & 1.00           & 1.50                  \\
       & 4   & 1522.91 & \multicolumn{1}{c}{-} & 363.96  & 16.38        & 248.18       & 1787.85 & 1.00           & 1.15                  \\ \hline
HFR    & 1   & 200.99  & 0.56                  & 2.71    & 3.14         & 25.45        & 230.15  & 0.94           & 0.72                  \\
       & 2   & 208.59  & 3.32                  & 12.50   & 3.74         & 46.32        & 261.97  & 1.10           & 1.58                  \\
       & 3   & 433.31  & 16.89                 & 54.76   & 8.17         & 111.03       & 569.40  & 1.27           & 1.91                  \\
       & 4   & 922.10  & 60.34                 & 195.94  & 16.71        & 248.13       & 1247.28 & 1.43           & 1.65                  \\ \hline
EFR    & 1   & 201.28  & 0.59                  & 2.61    & 3.19         & 25.86        & 230.93  & 0.94           & 0.72                  \\
       & 2   & 153.48  & 2.79                  & 12.32   & 4.50         & 53.55        & 214.32  & 1.35           & 1.93                  \\
       & 3   & 225.92  & 16.51                 & 52.08   & 8.67         & 114.95       & 366.05  & 1.98           & 2.97                  \\
       & 4   & 452.98  & 56.68                 & 189.40  & 16.52        & 247.93       & 774.11  & 2.31           & 2.66                  \\ \hline
\end{tabular}
\end{table}
% \usepackage{multirow}
Contours of Q-criterion are shown for simulations using $p=2$ and $p=4$ in Figure~\ref{fig:qcriterionturbcylside}, where the behaviour of these vortical structures can be observed. They result from instabilities caused by the complex shedding phenomena associated with this $\operatorname{Re}$. At higher polynomial degrees, finer turbulent structures are observed, which is expected due to the increased resolution of the $p=4$ method against $p=2$. Overall, contour results from both methods are in good agreement with each other and with the reference~\cite{zhaoDirectNumericalSimulation2009}. Results from this problem demonstrate the suitability of IMEX methods for problems of moderate stiffness.

\subsection{Turbulent Flow over a Multi-Element Airfoil}
\begin{figure}[htbp]
  \centering
    \begin{subfigure}[b]{0.98\textwidth}
    \includegraphics[width=\textwidth]{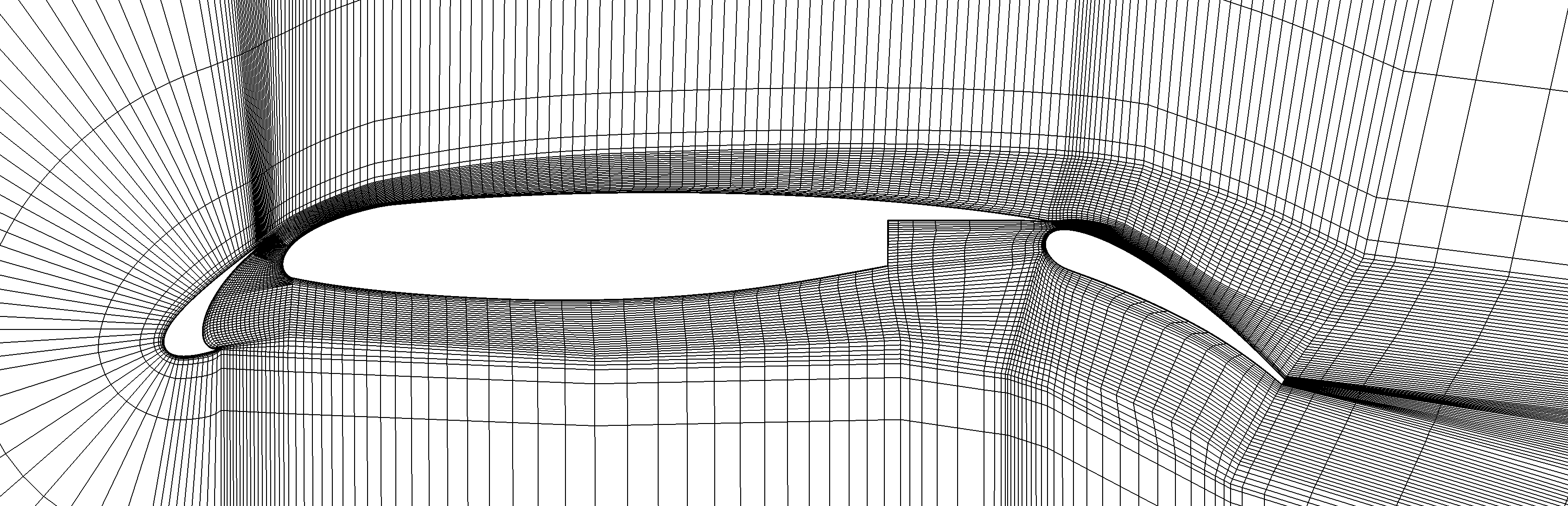}
    \caption{Multi-element view}
  \end{subfigure}
  \begin{subfigure}[b]{0.49\textwidth}
    \includegraphics[width=\textwidth]{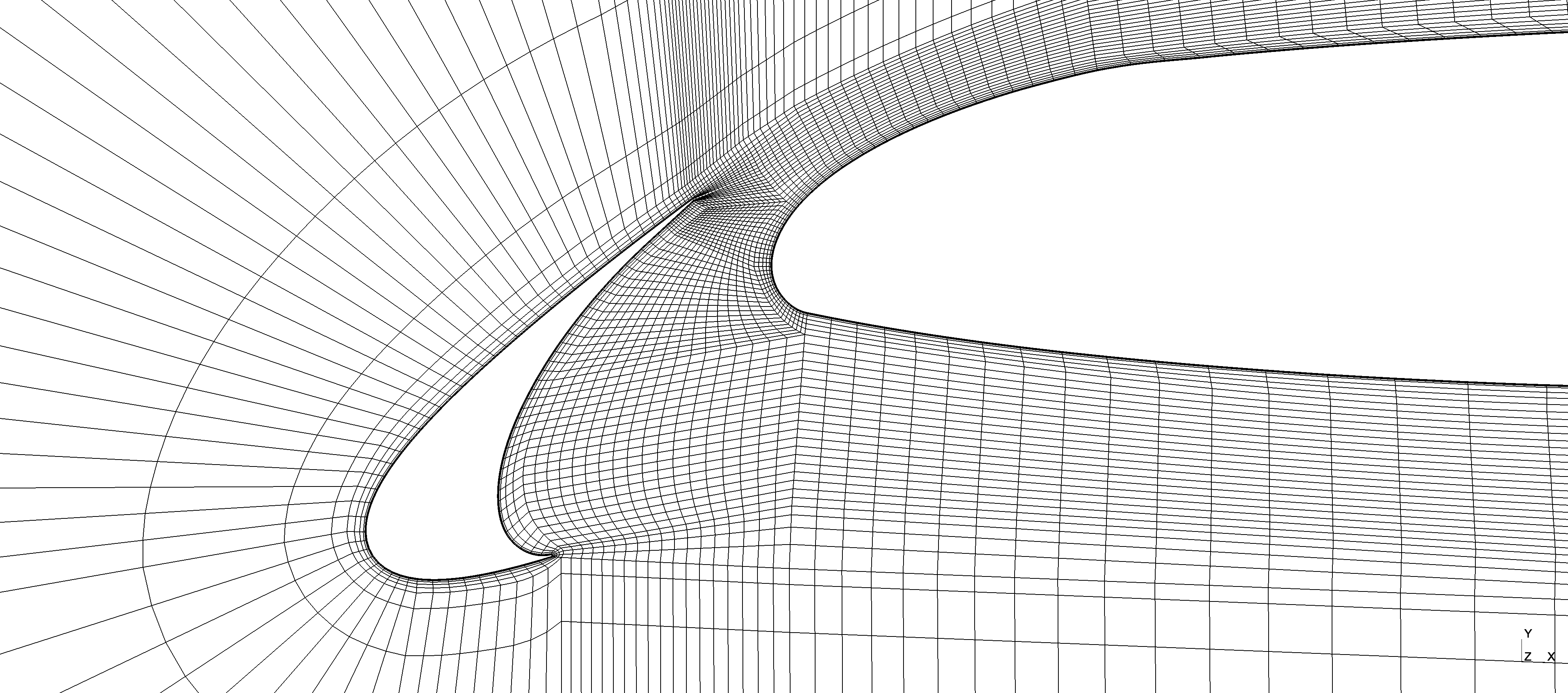}
    \caption{Slat view}
  \end{subfigure}
    \begin{subfigure}[b]{0.49\textwidth}
    \includegraphics[width=\textwidth]{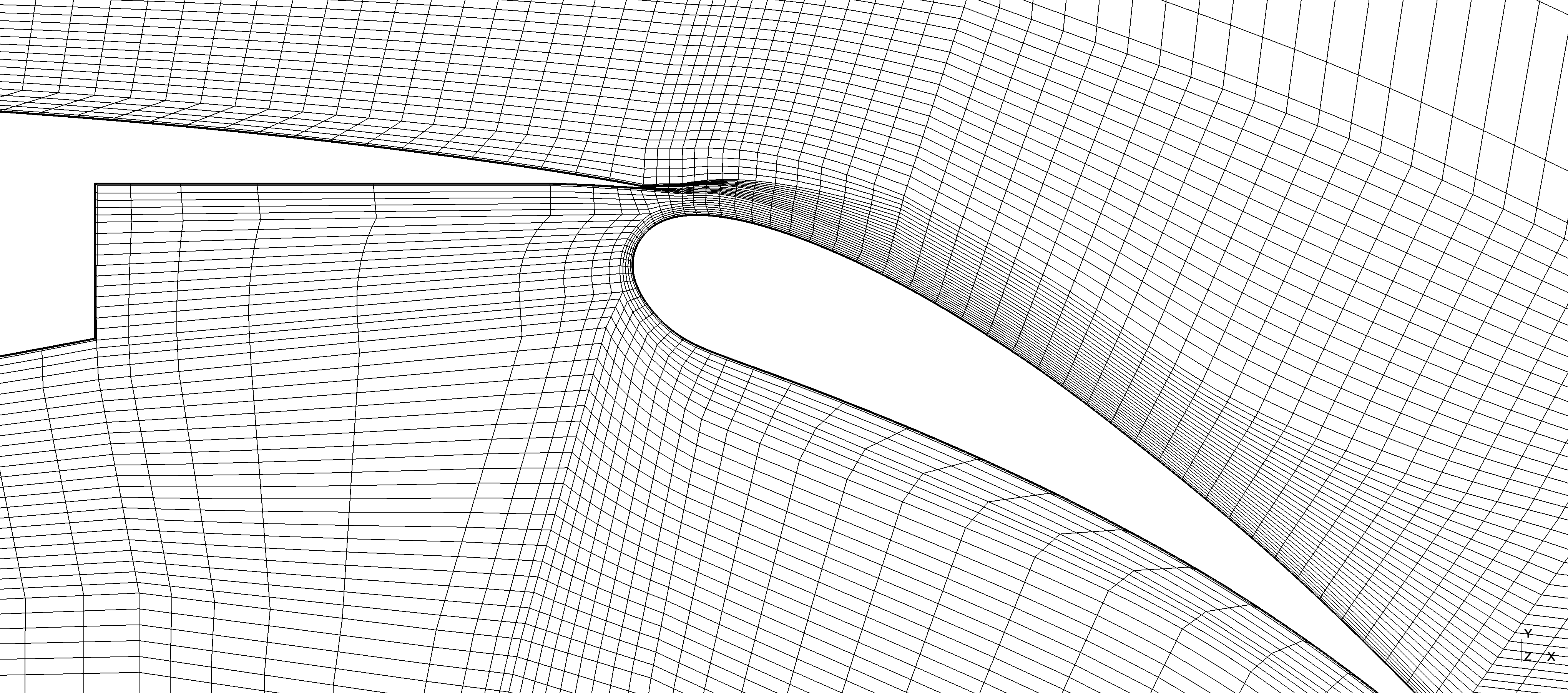}
    \caption{Flap view}
  \end{subfigure}
  \caption{Computational grid for the multi-element airfoil}
  \label{fig:30p30nmesh}
\end{figure}

Finally, this section presents wall-resolved large-eddy simulation of a multi-element 30P30N airfoil at $\operatorname{Re}=1.7\times10^6$ and an angle of attack 5.5 degrees. The reference chord is denoted $c$ and represents the stowed airfoil. This problem is commonly used in the aeroacoustics community, particularly within the AIAA Benchmark Problems in Airframe Noise Computations workshops~\cite{choudhari2015assessment}. This airfoil makes use of a slat at the leading edge to increase maximum lift, allowing the main wing to operate at higher angles of attack before stalling. These components introduce complex flow behaviour and significantly contribute to the acoustic field during the landing phase of an aircraft. Multiple researchers have experimentally~\cite{pascioni2016aeroacoustic} and numerically~\cite{ashton2016flow,terracol2014wall,choudhari2007effect} produced reference data, and an overview of results from the aforementioned workshop is available in~\cite{choudhari2015assessment}. This problem has been previously demonstrated in the context of wall-modeled LES and wall-resolved LES~\cite{choudhari2007effect,shi2018towards}, the latter with generally dense refinement within the slat cove region only and coarse everywhere else. The relatively high Reynolds number makes this problem computationally challenging. We generate a computational grid of 549280 hexahedral elements, shown in Figure~\ref{fig:30p30nmesh}, with zoomed-in views for the slat and the flap. The spanwise length is $c/9$, with 40 elements uniformly refined, which is slightly above the coarse grid with 30 layers used in~\cite{choudhari2007effect}. 

The entropically-damped artificial compressibility (EDAC) method is used in this problem, which resolved stability issues encountered with the compressible Navier-Stokes equations~\cite{BOLDUC2023105839}. This is an appropriate choice since the baseline Mach number for this problem is 0.1, which is within the incompressible range. We set the incompressibility factor to $\Theta=100$. This value was chosen to maintain a sensible time-step size and reduce the effects of artificial compressibility. See~\cite{trojak2022artificial} for a discussion of this parameter. We make use of a Rusanov-type Riemann solver, with stabilization computed from Davis estimates of the maximum eigenvalues of the EDAC equations~\cite{trojak2022artificial}. We include Appendix A with the resulting forms of the equations and parameters. For the hybridized formulation, we employ the following convective stabilization parameter
\begin{equation}
  s = \frac{3}{2}|\hat{\bm v}| + \hat d,\quad\quad \hat d^2 = \frac{\hat{\bm v}^2}{4} + \hat P + \Theta,
\end{equation}
which leads to an isotropic stabilization operator $s\bm I$. This stabilization leads to a new hybridized formulation of the EDAC equations, first used in this work. Simulations are run for 20 convective times, and the statistics are averaged for the last 10$t_c$. The time-step sizes considered here are $\Delta t/t_c=8\times10^{-5}$ for $p=1$, $\Delta t/t_c=3.125\times10^{-5}$ for $p=2$ and $\Delta t/t_c=1\times10^{-5}$ for $p=3$. Due to the underresolution caused by the coarseness of the grid, we applied regularization to the $p=3$ configuration via modal filtering. The choice of parameters in the modal filter was made heuristically by increasing the strength of the filter until stabilization was achieved~\cite{hamedi2022optimized}. We describe the filter and the parameters used in Appendix B.

\begin{figure}[htbp]
  \centering
  \includegraphics[width=0.6\textwidth]{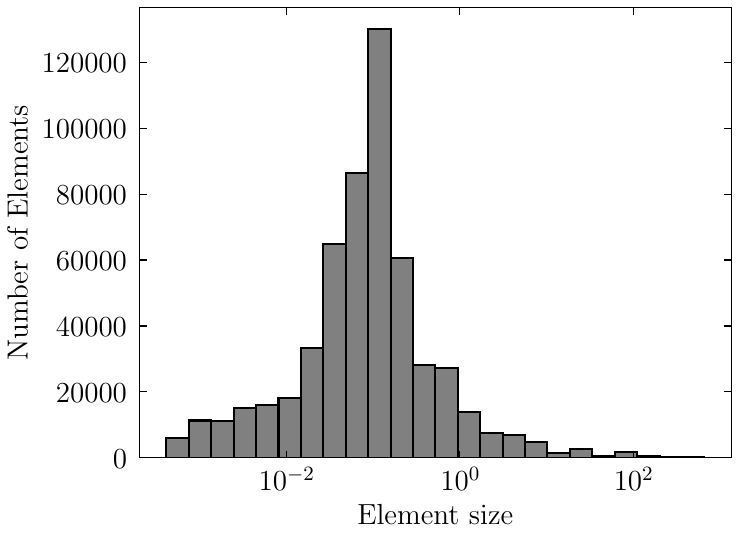}
  \caption{Distribution of element sizes for the multi-element airfoil grid}
  \label{fig:30p30nhistogram}
\end{figure}

\begin{figure}[htbp]
  \centering
    \begin{subfigure}[b]{0.49\textwidth}
    \includegraphics[width=\textwidth,trim={12cm 0cm 12cm 0cm},clip]{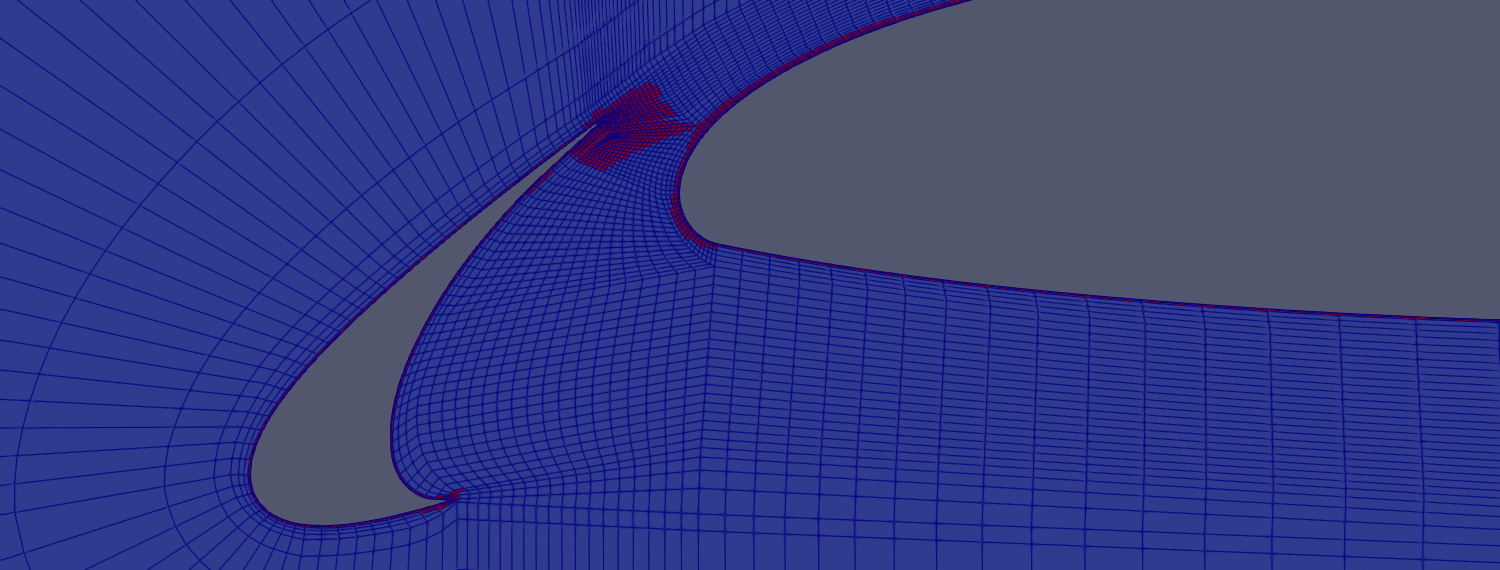}
    \caption{Slat cove zoom-in}
  \end{subfigure}
  \begin{subfigure}[b]{0.49\textwidth}
    \includegraphics[width=\textwidth,trim={12cm 0cm 12cm 0cm},clip]{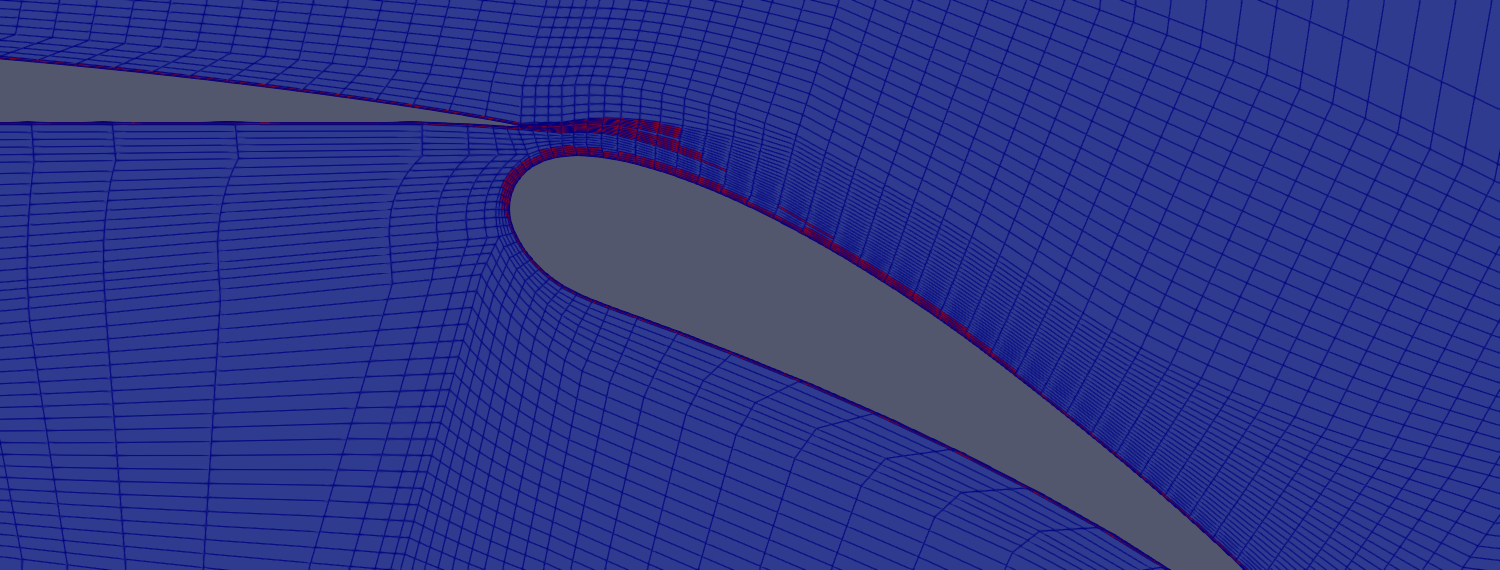}
    \caption{Flap zoom-in}
  \end{subfigure}
  \caption{Distribution of implicit and explicit elements after IMEX partitioning. Blue elements are solved explicitly, and red elements are solved implicitly}
  \label{fig:30p30nimexdomain}
\end{figure}

Following the promising results of the EFR method in the turbulent cylinder problem, we employ this method in the implicit portion of our IMEX domain but will also include performance results against the HFR and FR methods by computing a number of time steps in those. The implicit factor is chosen to solve $\sim$20\% of elements with the EFR method at $\operatorname{IF}=0.23$, consistent with our previous analyses in the cylinder problems. A histogram with the element size in the $x$-axis is provided in Figure~\ref{fig:30p30nhistogram}, showing element sizes with over four orders of magnitude in difference. The final distribution of implicit (red) and explicit (blue) elements in the domain partition is shown in Figure~\ref{fig:30p30nimexdomain}. Elements near the walls will be resolved implicitly, and medium to large elements will be solved explicitly. A tolerance of $10^{-4}$ is used to converge the unsteady implicit residuals.

\begin{figure}[htbp]
  \centering
  \begin{subfigure}[b]{0.49\textwidth}
    \includegraphics[width=\textwidth,trim={7cm 2cm 13cm 2cm},clip]{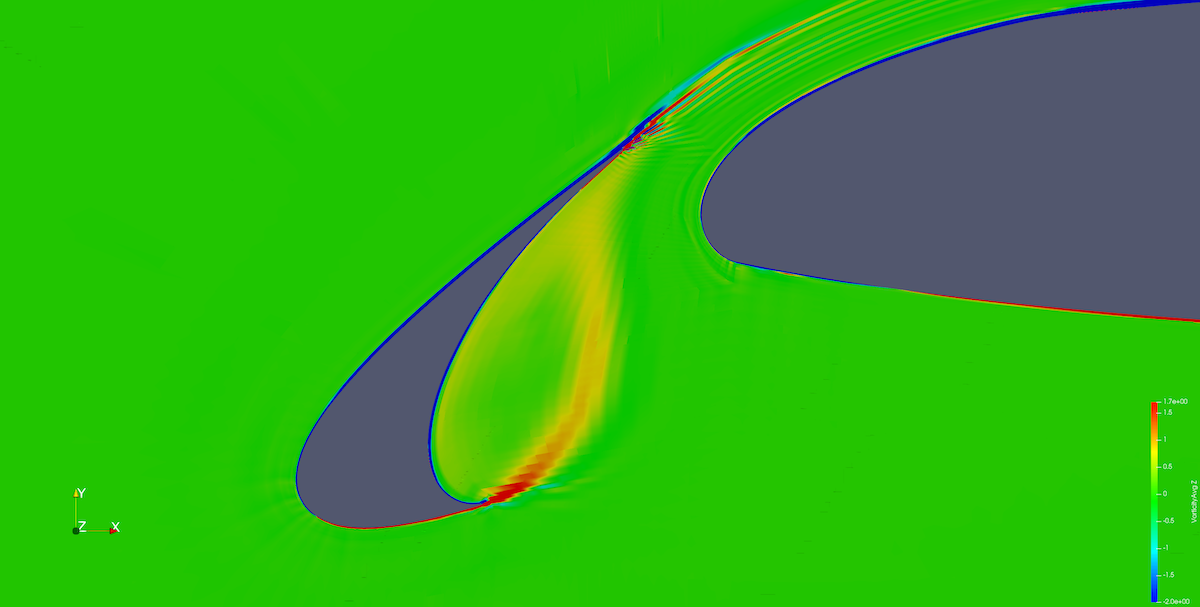}
    \caption{Averaged, $p=1$}
  \end{subfigure}
      \begin{subfigure}[b]{0.49\textwidth}
    \includegraphics[width=\textwidth,trim={7cm 2cm 13cm 2cm},clip]{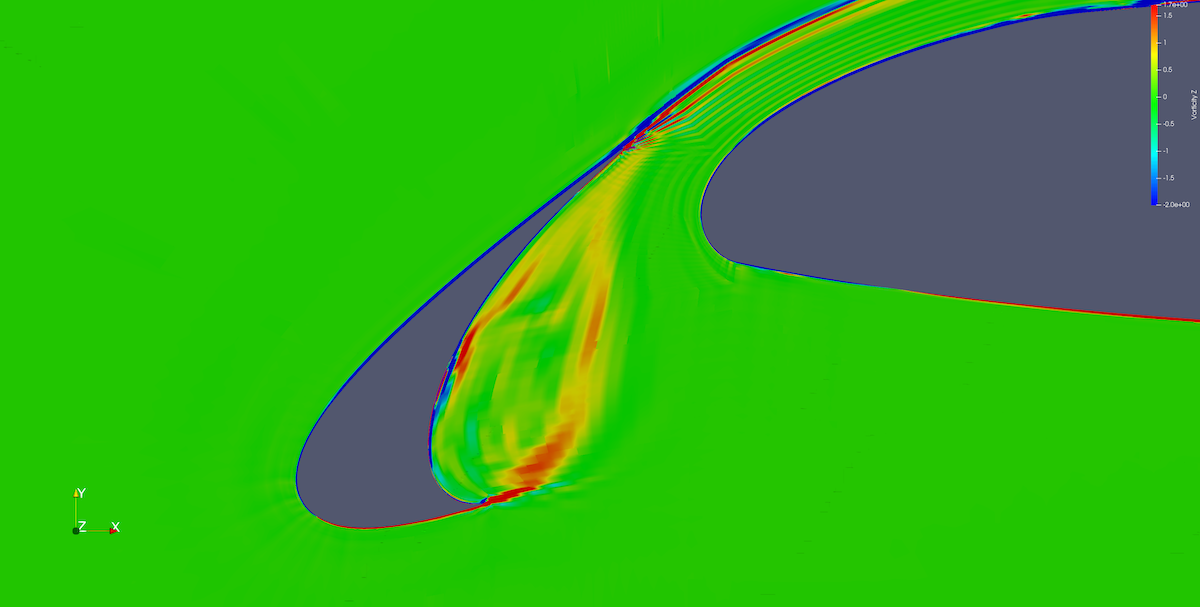}
    \caption{Instantaneous, $p=1$}
  \end{subfigure}

    \begin{subfigure}[b]{0.49\textwidth}
    \includegraphics[width=\textwidth,trim={7cm 2cm 13cm 2cm},clip]{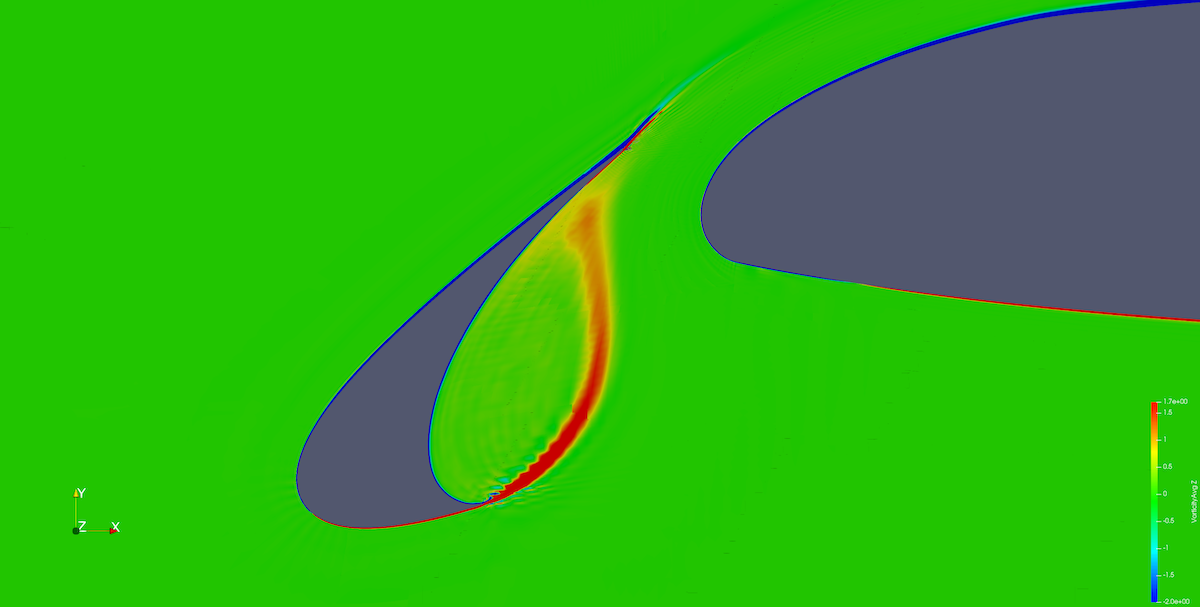}
    \caption{Averaged, $p=2$}
  \end{subfigure}
      \begin{subfigure}[b]{0.49\textwidth}
    \includegraphics[width=\textwidth,trim={7cm 2cm 13cm 2cm},clip]{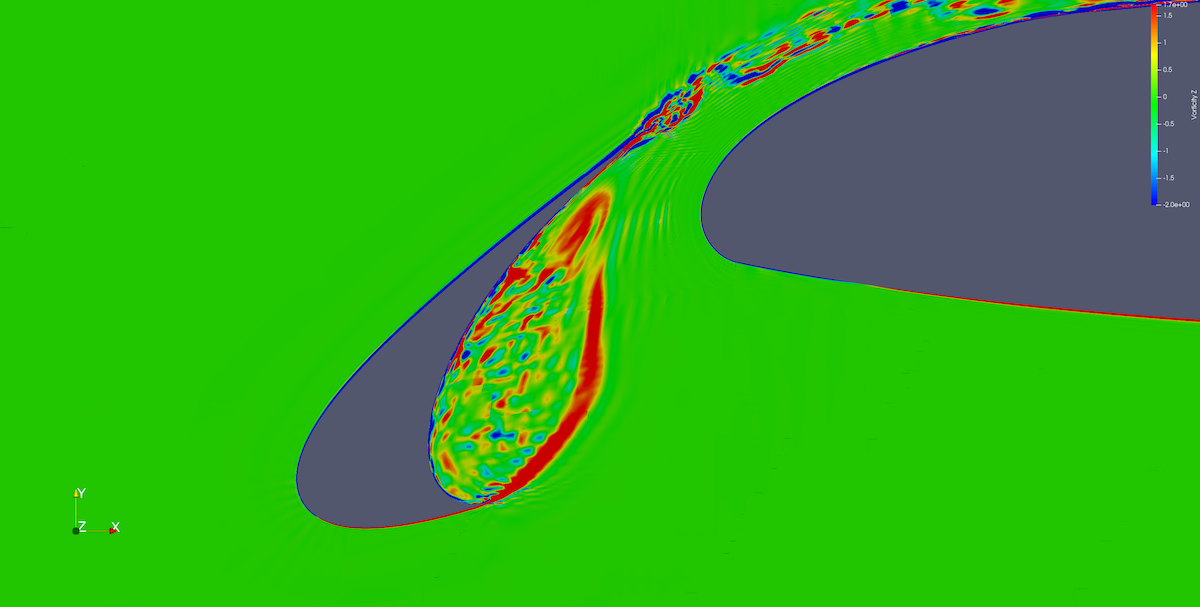}
    \caption{Instantaneous, $p=2$}
  \end{subfigure}
      \begin{subfigure}[b]{0.49\textwidth}
    \includegraphics[width=\textwidth,trim={7cm 2cm 13cm 2cm},clip]{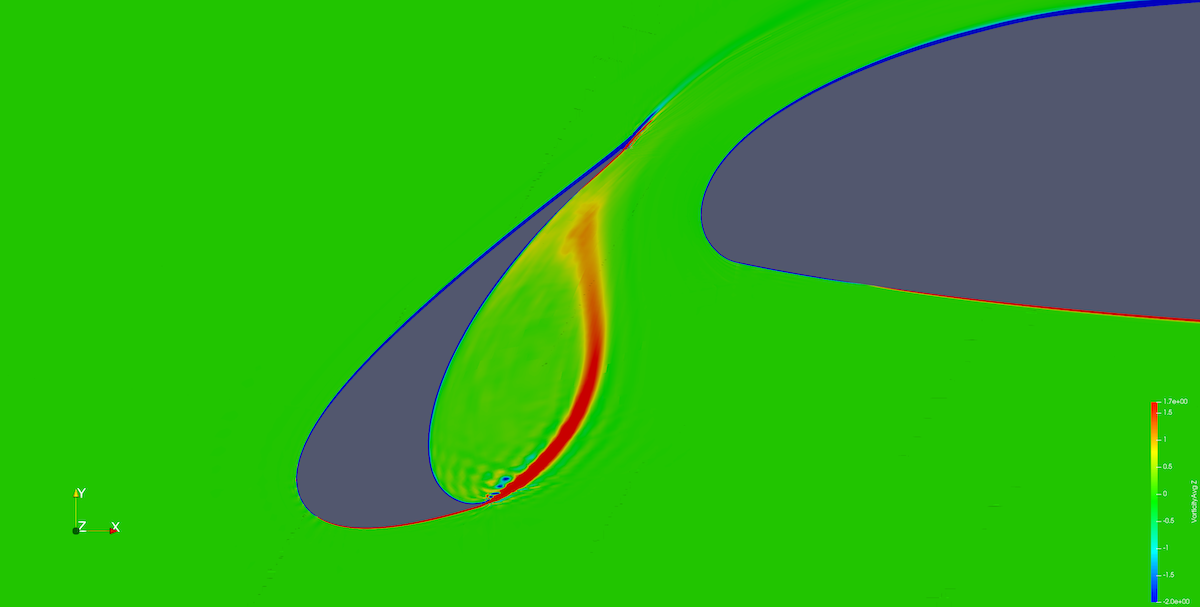} % modify
    \caption{Averaged, $p=3$}
  \end{subfigure}
      \begin{subfigure}[b]{0.49\textwidth}
    \includegraphics[width=\textwidth,trim={7cm 2cm 13cm 2cm},clip]{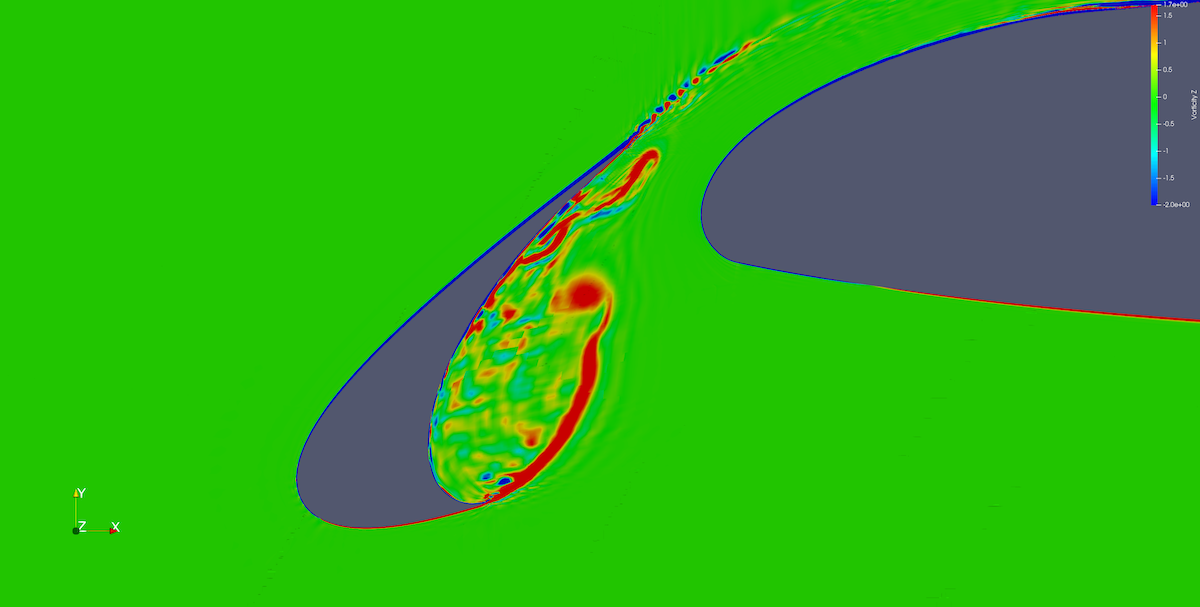} % modify
    \caption{Instantaneous, $p=3$}
  \end{subfigure}
  \caption{Vorticity contours for the multi-element airfoil problem}
  \label{fig:30p30nvorticity}
\end{figure}

Contours of averaged vorticity are shown in Figure~\ref{fig:30p30nvorticity} for $p=1$ to $p=3$ simulations. Results from the $p=1$ simulations are highly dissipative and display smeared-out regions of vorticity. Simulations at $p=2$ and $p=3$ are already in good agreement with the PIV visualizations from Pascioni et al.~\cite{pascioni2014experimental}, where the expected detached shear layer emerging from the slat is observed with increased definition. In the instantaneous plots, similar behaviour is observed between the $p=1$ and $p=2$ simulations, resulting in the latter being a more accurate representation of the complex vortex interaction within the slat cove. The instantaneous snapshots for $p=3$ display larger structures, which can be due to the dissipation caused by the strength of the filter.

\begin{figure}[htbp]
  \centering
    \begin{subfigure}[b]{0.7\textwidth}
    \includegraphics[width=\textwidth,trim={5cm 2cm 5cm 2cm},clip]{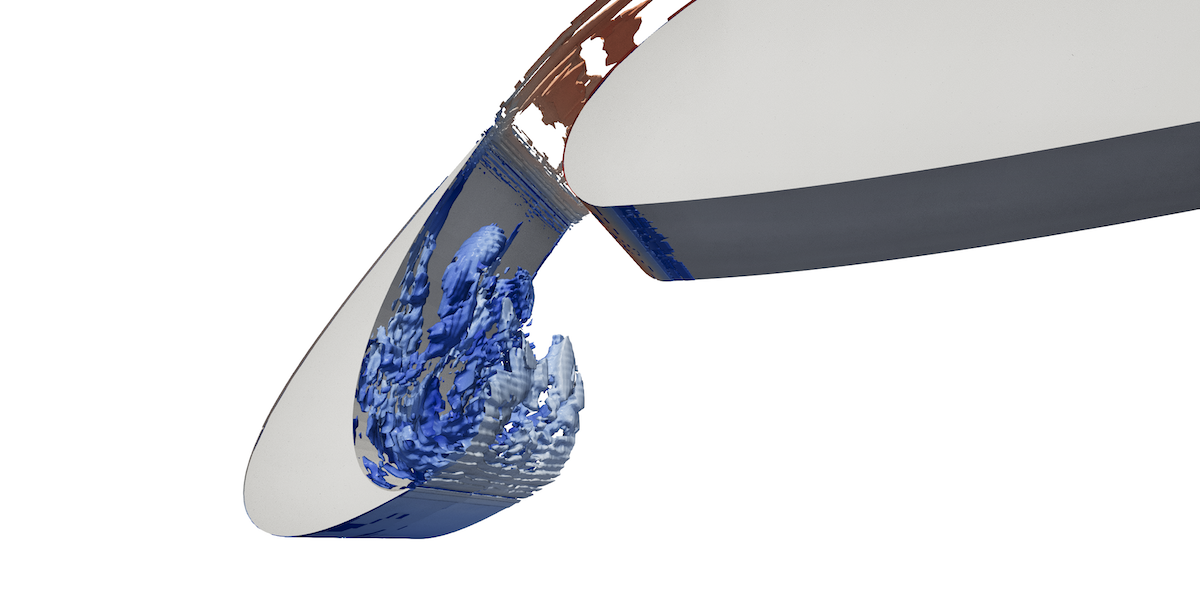}
    \caption{$p=1$}
  \end{subfigure}
  \begin{subfigure}[b]{0.7\textwidth}
    \includegraphics[width=\textwidth,trim={5cm 2cm 5cm 2cm},clip]{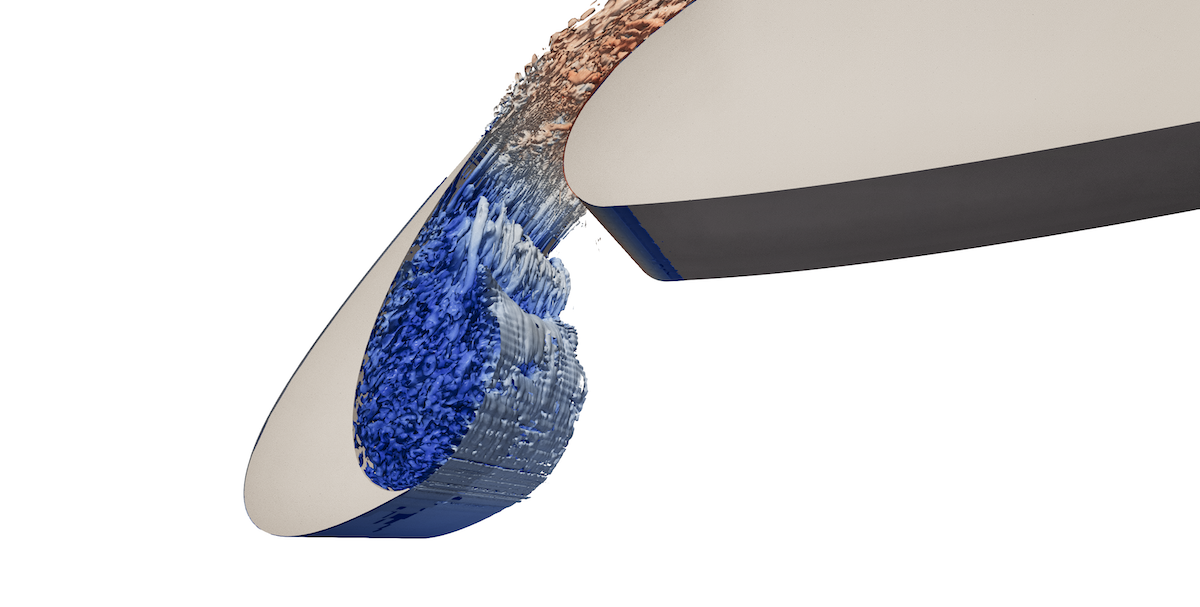}
    \caption{$p=2$}
  \end{subfigure}
    \begin{subfigure}[b]{0.7\textwidth}
    \includegraphics[width=\textwidth,trim={5cm 2cm 5cm 2cm},clip]{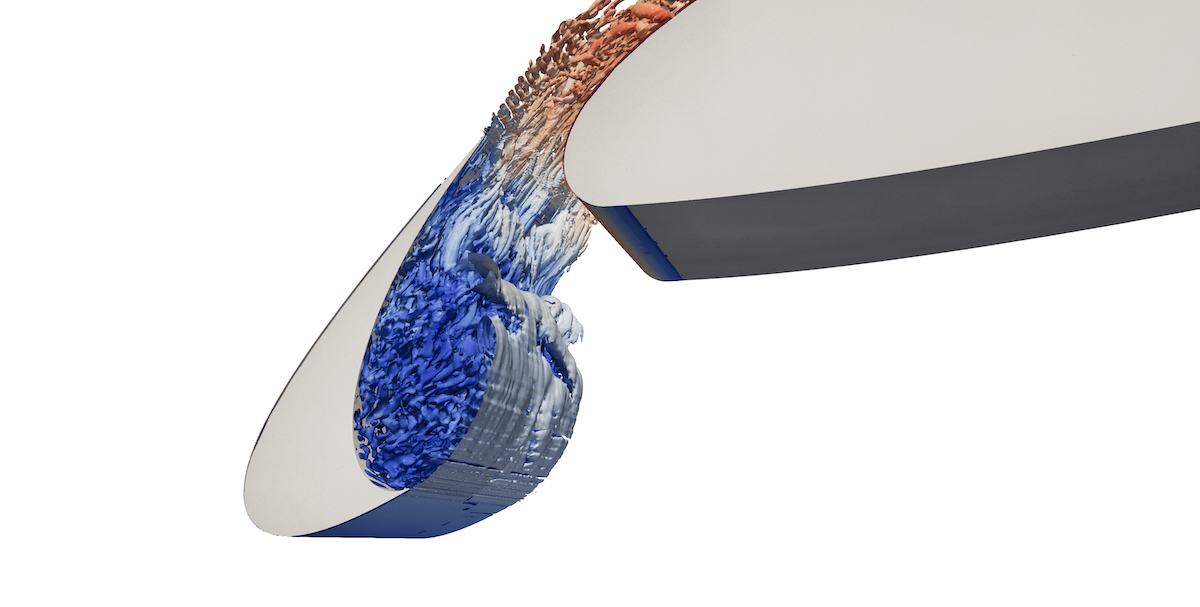}
    \caption{$p=3$}
  \end{subfigure}
  \caption{Q-criterion contours zoomed in the slat for the multi-element airfoil coloured by streamwise vorticity}
  \label{fig:qcriterion30p30nslat}
\end{figure}

A similar conclusion can be drawn from the contours of Q-criterion, which are shown for a zoomed-in portion within the slat cove in Figure~\ref{fig:qcriterion30p30nslat}. The second-order simulations resolve a smaller range of scales compared to third and fourth orders. At the chosen value of Q-criterion, many more structures are shown in the higher-order simulation emanating from the gap between the cove and airfoil all the way to the upper side of the flap downstream. For $p=1$, the detail in the turbulent structures is small, whereas much more detail can be observed for the higher-order results.

\begin{figure}[htbp]
  \centering
   \includegraphics[width=0.65\textwidth]{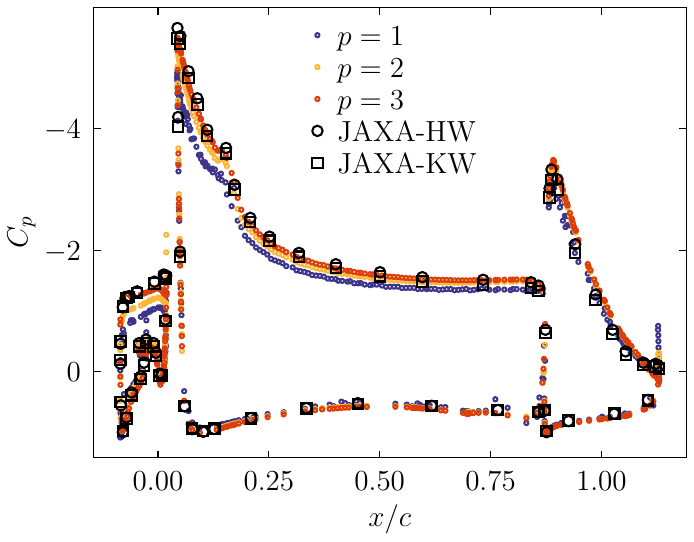}
  \caption{Pressure coefficient $C_p$ for the multi-element airfoil problem}
  \label{fig:30p30ncp}
\end{figure}

\begin{figure}[htbp]
  \centering
   \includegraphics[width=0.65\textwidth]{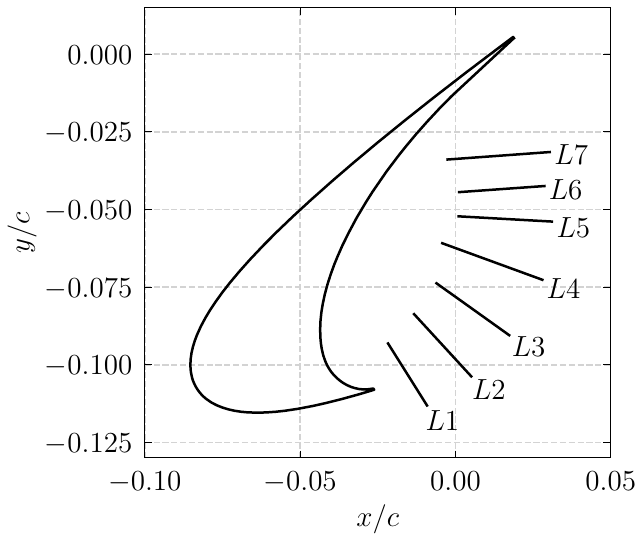}
  \caption{Line plot locations in the vicinity of the slat}
  \label{fig:linecuts}
\end{figure}

\begin{figure}[htbp]
  \centering
  \begin{subfigure}[b]{0.49\textwidth}
    \includegraphics[width=\textwidth]{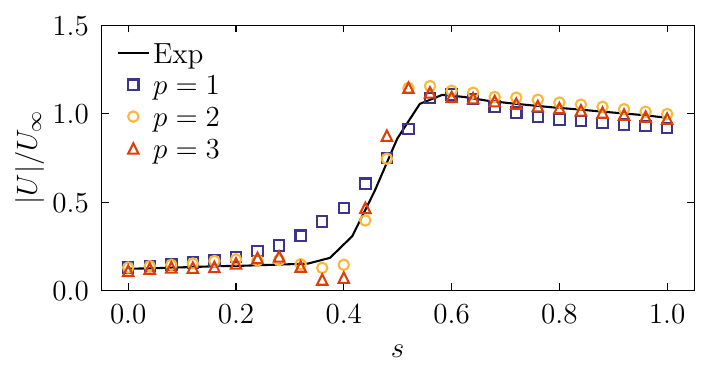}
    \caption{L1}
  \end{subfigure}
    \begin{subfigure}[b]{0.49\textwidth}
    \includegraphics[width=\textwidth]{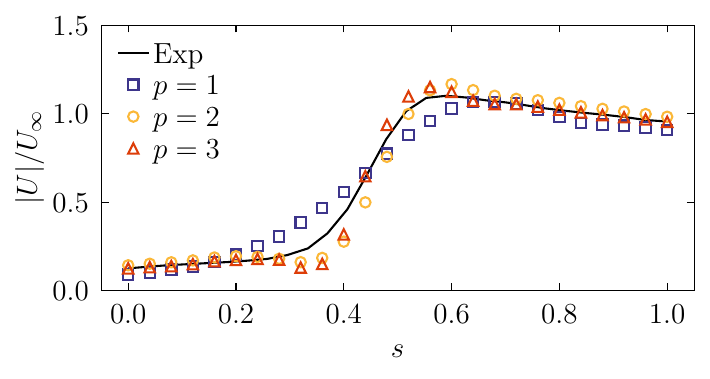}
    \caption{L2}
  \end{subfigure}
    \begin{subfigure}[b]{0.49\textwidth}
    \includegraphics[width=\textwidth]{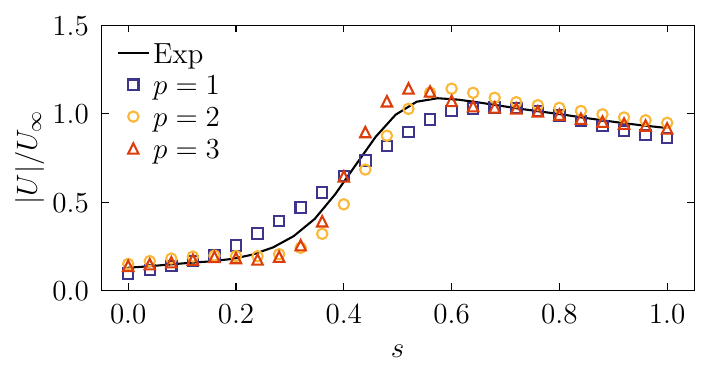}
    \caption{L3}
  \end{subfigure}
    \begin{subfigure}[b]{0.49\textwidth}
    \includegraphics[width=\textwidth]{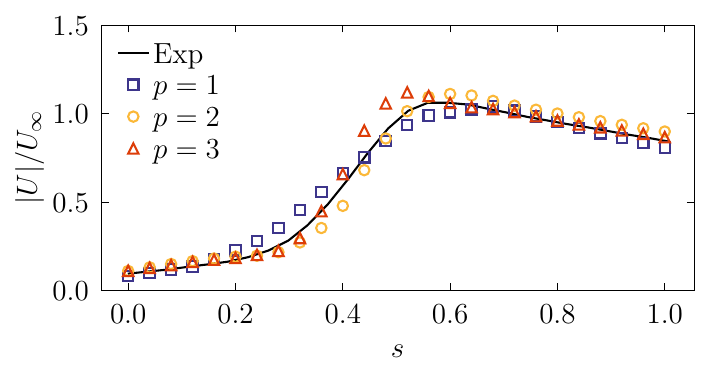}
    \caption{L4}
  \end{subfigure}
    \begin{subfigure}[b]{0.49\textwidth}
    \includegraphics[width=\textwidth]{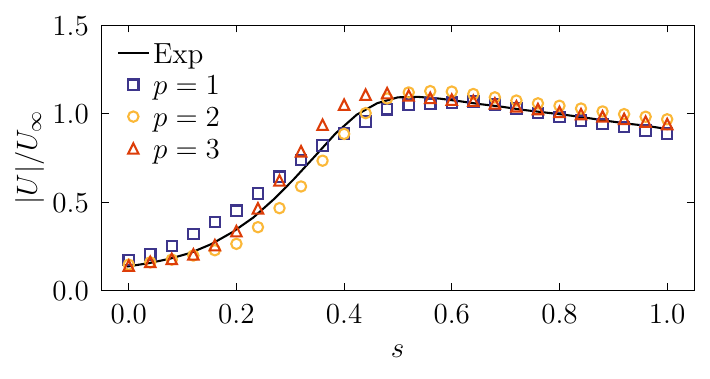}
    \caption{L5}
  \end{subfigure}
    \begin{subfigure}[b]{0.49\textwidth}
    \includegraphics[width=\textwidth]{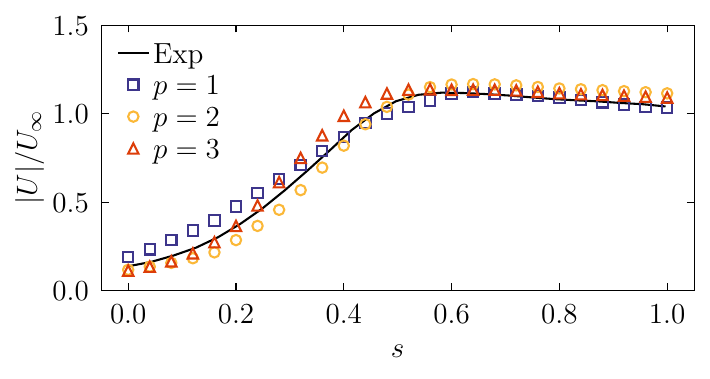}
    \caption{L6}
  \end{subfigure}
    \begin{subfigure}[b]{0.49\textwidth}
    \includegraphics[width=\textwidth]{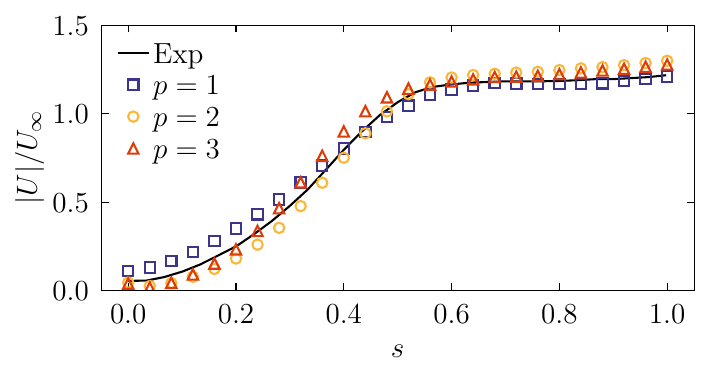}
    \caption{L7}
  \end{subfigure}
  \caption{Plots of sampled normalized velocity magnitude for the multi-element airfoil across slat shear layer}
  \label{fig:30p30nlines}
\end{figure}

Time-averaged plots of the pressure coefficient are shown in Figure~\ref{fig:30p30ncp}. Results are compared against experimental data by Florida State University~\cite{pascioni2016aeroacoustic} and from Muyarama et al. ~\cite{murayama2014experimental}. We note that these experimental values were originally performed in wind tunnel facilities at different angles of attack. They were carried out in closed-wall wind tunnels with significant end-wall effects, altering the effective attack angle. However, they were later compared to numerical simulations and were deemed appropriate as reference data at an angle of attack $5.5$ degrees~\cite{pascioni2016aeroacoustic}. Increasing the order to $p=2$ improves the agreement between the current results and the reference data, particularly for the main airfoil and the flap. The $C_p$ plots are close to the reference data for $p=2$ in these two regions, but the slat remains still quite underpredicted. Results for the $p=3$ simulation are closer to the reference data for the slat. We also compare velocity profiles along lines normal to the shear layer, shown in Figure~\ref{fig:linecuts}. Results for the three computations considered in this section are shown in Figure~\ref{fig:30p30nlines}, namely $p=1$, $p=2$ and $p=3$. Overall, an improved agreement can be observed for the third-order results, especially in proximity to the shear layer. In this region, the higher-order results follow the increases in velocity to the freestream conditions. The low-order method shows more dissipated behaviour in these regions, which is expected and consistent with the numerical error of these schemes. The $p=2$ results show oscillatory behaviour for the higher-order results at the crossing of the shear layer in L1. In general, there is still some discrepancy in areas close to the mean quantities, especially for $p=1$.

Finally, we measure computational performance. Specifically, the time spent in 100 time steps for FR, HFR, and EFR methods is computed. These timing metrics are shown and compared with the standard FR IMEX approach and with their explicit counterparts when the implicit factor is zero, which results in an optimal explicit-Runge Kutta method. Results are shown in Table~\ref{tab:30p30nperformanceimex}. These runtimes were computed on 2.4GHz AMD Rome 7532 CPUS using 1024 cores. Similar to previous results in the cylinder, significant time is spent on the computation of the global problems, especially at higher order, accounting for 80\% of the time in FR, over 70\% in HFR, and over 50\% in EFR, with most of the remaining time in the explicit portion. Compared to the standard FR IMEX approach, EFR IMEX methods achieved 1.5 to 4.1 times faster simulations, representing at least 15 times faster than an explicit FR simulation of this problem. Hence, we have demonstrated that hybridized IMEX methods can significantly speed up computations at high Reynolds numbers. 
\begin{table}[htbp]
\centering
\caption{Summary of performance metrics for the multi-element airfoil problem for 100 time steps}
\label{tab:30p30nperformanceimex}
\begin{tabular}{lrrrrrrrrr}
\hline
Scheme & $p$ & $t_{G}$ & $t_{L}$               & $t_{J}$ & $t_{R^{\text{im}}}$ & $t_{R^{\text{ex}}}$ & $t_w$   & $t_w/t_{w}^{\text{FR}}$ & $t_w/t_{w}^{\text{ex}}$ \\ \hline
FR     & 1   & 14.66   & \multicolumn{1}{c}{-} & 4.72    & 0.65         & 7.79         & 27.81   & 1.00           & 35.83                 \\
       & 2   & 101.41  & \multicolumn{1}{c}{-} & 19.48   & 2.04         & 28.04        & 150.97  & 1.00           & 24.56                 \\
       & 3   & 983.59  & \multicolumn{1}{c}{-} & 134.93  & 6.28         & 55.52        & 1180.32 & 1.00           & 4.03                  \\ \hline
HFR    & 1   & 54.61   & 1.08                  & 1.29    & 0.63         & 7.51         & 65.13   & 0.43           & 15.33                 \\
       & 2   & 218.56  & 6.46                  & 8.21    & 2.33         & 28.49        & 264.04  & 0.57           & 14.06                 \\
       & 3   & 370.49  & 15.90                 & 60.97   & 6.42         & 52.03        & 505.81  & 2.33           & 9.40                  \\ \hline
EFR    & 1   & 7.72    & 1.25                  & 1.13    & 0.67         & 7.79         & 18.56   & 1.50           & 53.80                 \\
       & 2   & 62.86   & 6.50                  & 7.50    & 2.40         & 27.97        & 107.21  & 1.41           & 34.62                 \\
       & 3   & 153.08  & 15.61                 & 59.02   & 6.27         & 51.18        & 285.17  & 4.14           & 16.67                 \\ \hline
\end{tabular}
\end{table}

\section{Conclusions}\label{sec:conclusion}
Novel hybridized implicit-explicit methods were developed for large-scale simulations with geometry-induced stiffness. Combining hybridized and standard FR formulations can be done conservatively using the FR formulation's fluxes at the interface. Via a performance analysis in two and three dimensions, it was shown that for problems with moderate stiffness, obtained speedups are significant in two dimensions. For a laminar cylinder problem, speedups over 6$\times$ were obtained at the highest polynomial degrees, comparing EFR against FR IMEX methods. In the three-dimensional moderately stiff setting via a turbulent cylinder at $\operatorname{Re}=1,000$, these speedups are only observed at higher polynomial degrees $p>3$ for HFR and $p\geq 2$ for EFR. However, results for a multi-element airfoil at $\operatorname{Re}=1.7\times10^6$ demonstrated that in problems where geometry-induced stiffness is a significant contributor, such as high Reynolds numbers, the speedups of IMEX methods are significant. Performance speedups of EFR methods against FR-IMEX schemes were in excess of four, resulting in simulations at least 15 times faster than explicit counterpart formulations. IMEX with hybridization demonstrated significant potential in flow computations using the FR approach. The properties of this method can be expanded to simulate additional applications that benefit from domain subdivisions with localized stiffness. This can be, for instance, fluid-solid interaction problems, where hybridization can be used in the solid portion and FR in the fluid portion. This will allow a seamless extension of flow solvers to more complex applications.

\section*{Acknowledgements}
\noindent We acknowledge the support of the Natural Sciences and Engineering Research Council of Canada (NSERC), [RGPAS-2017-507988, RGPIN-2017-06773], Concordia University via the Team Seed pro- gram and the Fonds de Recherche du Quebec - Nature et Technologie (FRQNT) via a B2X scholarship. This research was enabled in part by support provided the Digital Research Alliance of Canada (www.alliancecan.ca) via a Resources for Research Groups allocation.

\section*{Data Statement}
\noindent Data relating to the results in this manuscript can be downloaded from the publication's website under a CC BY-NC-ND 4.0 license.

\pagebreak

\bibliography{manuscript}

\appendix
\section{EDAC Equations}
As the Mach number decreases, the disparity between the entropy and the acoustic wave speeds in the flow becomes significant. Hence, for these problems, solutions via compressible Navier-Stokes become quite challenging. Artificial compressibility methods (ACM)~\cite{chorin1997numerical} provide mechanisms toward the divergence-free condition of the velocity field via pseudo-time~\cite{jameson1991time}, which reduces the stiffness of the pressure and velocity field toward a divergence-free result. These methods enable explicit time stepping, contrary to the typical solutions of the incompressible Navier-Stokes equations~\cite{panton2013incompressible}, which requires solving a Poisson problem to obtain the pressure field.

The entropically-damped artificial compressibility (EDAC) method of Clausen~\cite{clausen2013entropically} achieves closure by minimizing density variations rather than setting a constant entropy constraint, such as in the ACM method. The resulting equation has a pressure diffusion operator. The pressure evolution is hence dictated by
\begin{equation}
\frac{\partial P}{\partial t} + \bm v \cdot \nabla P  + \frac{1}{\operatorname{M}^2} \nabla\cdot \bm v - \frac{1}{\operatorname{Re}} \nabla^2 P = 0.
\end{equation}
We can describe a general form of these equations in the form of a convection-diffusion equation 
\begin{equation}
    \frac{\partial \bm{u}}{\partial t}
    + \sum_{i=1}^{d} \frac{\partial}{\partial x_i} \left[\bm{F}^{(c)}_{i}(\bm u)-\bm{F}^{(v)}_{i}(\bm u, \nabla \bm u)\right]
   = \bm{s}\quad\text{in } \Omega,
\end{equation}
Hence, the vector of conserved variables for the EDAC equations is given by
\begin{equation}
  \bm u = 
  \begin{bmatrix}
   P \\ v_x \\ v_y \\ v_z 
  \end{bmatrix},
\end{equation}
as well as convective and viscous flux vectors defined by
\begin{equation}
   \bm F_i^{(c)} = 
  \begin{bmatrix}
    v_i(P + \Theta)\\
    v_i v_x + P\delta_{ix} \\
    v_i v_y + P\delta_{iy} \\
    v_i v_z + P\delta_{iz} \\
  \end{bmatrix},
  \quad
  \bm{F}^{(v)}_{i} = \nu
  \frac{\partial}{\partial x_i}
  \begin{bmatrix}
    P \\
    v_x \\
    v_z \\
    v_y \\
  \end{bmatrix},
\end{equation}
where $\Theta=1/\operatorname{M}^2$ with $\operatorname{M}$ an artificial Mach number, and $\nu=\frac{1}{\operatorname{Re}}$. These equations have been successfully applied in the context of FR methods~\cite{trojak2022artificial}, showing that they are an effective alternative to the ACM method, providing faster results and simpler implementation. Higher values of $\Theta$ introduce stiffness in the problem, but provide a more accurate approximation of the divergence-free condition. The maximum stable time step size is also influenced by this parameter in explicit numerical schemes. We make use of the EDAC equations in the simulation of a multi-element airfoil problem, where we take advantage of their explicit form to employ hybridized implicit-explicit methods.

\section{Modal Filtering}
The filtering approach used in the multi-element airfoil problem makes use of an operator of the form
\begin{equation}
  \bm F = \bm V \bm \Lambda^* \bm V^{-1},
\end{equation}
where $\bm V$ is a vandermonde matrix $V_{ij}=\Psi_{j-1}(\tilde{\bm x}_i)$ with $\Psi_j$ an orthonormal Legendre polynomial of degree $j$. $\bm \Lambda^*$ is the modal filtering matrix, normalized to be independent of the time-step size~\cite{hamedi2022optimized}. The filtering matrix before normalization is a diagonal matrix with entries
\begin{equation}
\bm \Lambda_{ii} = \sigma(\eta),
\end{equation}
where $\eta$ is the sum of the exponents in the leading term of the corresponding orthonormal basis function and $\sigma$ is the filtering function defined by
\begin{equation}
  \sigma(\eta) = 
  \begin{cases}
    1& 0 \leq \eta < \eta_c, \\
    \exp\left[-\alpha\left(\frac{\eta-\eta_c}{\eta_{\max}-\eta_c}\right)^s\right]& \eta_c\leq \eta \leq \eta_{\max}, \\
    0&\eta>\eta_{\max}.
  \end{cases}
\end{equation}
In this function, $\alpha$, $s$ are the damping and strength parameters. $\eta_{\max}$ is the maximum exponent in the orthonormal basis, and $\eta_c$ is the cut-off degree. For the simulation of the multi-element airfoil, we achieved stability after setting $\alpha=100$, $s=1$, and $\eta_c=3$.

\end{document}